\documentclass[12pt]{article}
\usepackage{fullpage,latexsym,amsfonts,amsthm,amsmath,amscd,amssymb,color}
\usepackage[dvips]{graphicx}
\usepackage{epic}

\newtheorem{lemma}{Lemma}[section]
\newtheorem{thm}[lemma]{Theorem}
\newtheorem{rem}[lemma]{Remark}

\newtheorem{example}[lemma]{Example}
\newtheorem{defn}[lemma]{Definition}

\newcommand{\ad}{\operatorname{ad}}
\newcommand{\Ad}{\operatorname{Ad}}

\def\C{{\mathbb C}}
\def\R{{\mathbb R}}

\def\E{{\cal E}}

\def\su{{\mathfrak{su}}}

\def\u{{\mathfrak{u}}}
\def\a{{\mathfrak{a}}}
\def\n{{\mathfrak{n}}}

\def\b{{\mathfrak{b}}}

\def\k{{\bf k}}

\def\g{\mathfrak{g}}

\def\h{{\mathfrak h}}

\def\t{\tilde}
\def\t{\mathfrak{t}}

\def\u{{\mathfrak u}}
\def\su{{\mathfrak {su}}}

\def\Omegakks{{\Omega_{kks}}}
\def\idempotent{{\epsilon}}

\newcommand{\BB}{\mathbb}

\newcommand{\bea}{\begin{eqnarray}}
\newcommand{\eea}{\end{eqnarray}}
\newcommand{\nn}{\nonumber}
\newcommand{\Tr}{\textrm{Tr}}
\newcommand{\bra}{\langle}
\newcommand{\ket}{\rangle}

\def\Im{{\rm Im}}
\def\prg{{\rm pr_{\u}}}
\def\prgsplus{{\rm pr_{\b_+}}}
\def\prgsminus{{\rm pr_{\b_-}}}
\def\divergence{{\rm div}}
\def\E{{\cal E}}

\def\Im{{\rm Im}}
\def\prg{{\rm pr_{\g}}}
\def\prgsplus{{\rm pr_{\b_+}}}
\def\prgsminus{{\rm pr_{\b_-}}}
\def\divergence{{\rm div}}
\def\sp{{\mathfrak {sp}}}
\def\so{{\mathfrak {so}}}

\DeclareMathAlphabet{\mathpzc}{OT1}{pzc}{m}{it}

\newtheorem{theorem}{Theorem}[section]
\newtheorem{proposition}[theorem]{Proposition}

\begin{document}

\title{\begin{flushright} \small
UUITP-08/15
 \end{flushright}
\bigskip 
\bigskip
Complete integrability from Poisson-Nijenhuis structures on compact hermitian symmetric spaces}

\author{F. Bonechi\footnote{\small INFN Sezione di Firenze, email: francesco.bonechi@fi.infn.it} ,
J. Qiu \footnote{\small Department of Mathematics, Uppsala University, Max-Planck-Institut f\"ur Mathematik, email: jian.qiu@math.uu.se},
~and
M. Tarlini\footnote{\small INFN Sezione di Firenze,  email: marco.tarlini@fi.infn.it}
}

\maketitle

\begin{abstract}
We study a class of Poisson-Nijenhuis systems defined on compact hermitian symmetric spaces, where the Nijenhuis tensor is defined as the
composition of Kirillov-Konstant-Souriau symplectic form with the so called Bruhat-Poisson structure. We determine its spectrum. In the case 
of Grassmannians the eigenvalues are the Gelfand-Tsetlin variables. We introduce the abelian algebra of collective 
hamiltonians defined by a chain of nested subalgebras and prove complete integrability. By construction, these models are integrable with respect to both Poisson structures. 
The eigenvalues of the Nijenhuis tensor are a choice of 
action variables. Our proof relies on an explicit formula for the contravariant connection defined on vector bundles that are Poisson with respect 
to the Bruhat-Poisson structure.
\end{abstract}

\thispagestyle{empty}

\eject
\section{Introduction}
Flag manifolds can be considered as homogeneous spaces of compact matrix group; when considered as coadjoint orbits they are endowed with
the {\it Kirillov-Konstant-Souriau} symplectic form $\Omegakks$. By fixing the standard Poisson-Lie structure on the matrix group, the quotient map
induces the {\it Bruhat-Poisson} structure $\pi_0$. It was shown in \cite{KRR} that the two Poisson structures, the inverse of the KKS symplectic form
and the Bruhat-Poisson, are compatible, {\it i.e.} their Schouten bracket vanishes, if and only if the flag manifold is a compact hermitian symmetric space.
This fact implies that there exists a {\it Poisson-Nijenhuis}
structure and most importantly there exists an integrable model admitting a {\it bihamiltonian} description.
In this paper we compute the eigenvalues of the Nijenhuis operator $N=\pi_0\circ\Omegakks$ for the cases of classical groups; these eigenvalues give
a specific choice of action variables.

In
\cite{KRR} and in \cite{Foth} it was shown that for complex projective spaces
these eigenvalues are given by the hamiltonians corresponding to fixing a certain basis of the torus; in particular it was noticed that they are
actually the Gelfand-Tsetlin variables. Moreover it was announced but not proved that this is true for all Grassmannians. This paper aims to fill this gap and 
generalize to the other cases.

Our motivation for understanding the properties of this integrable model comes from the problem of quantizing the symplectic groupoid integrating
the Bruhat-Poisson structure. This project was started in
\cite{BCST2} for $\C P_1$ and developed in \cite{BCQT} for $\C P_n$. The main idea is that thanks
to the groupoid structure we can use polarizations of the symplectic groupoid that are quite singular from the point of view of geometric quantization:
indeed we can consider real polarizations that induce on the space of lagrangian leaves the structure of topological groupoid. This is enough for defining the convolution
algebra from the groupoid of Bohr-Sommerfeld leaves (provided it admits a Haar system, which is true if,
for instance, it is \'etale). This observation led us in \cite{BCQT} to introduce the notion
{\it multiplicative integrability of the modular function}. The modular function is the groupoid cocycle that integrates the modular vector field of the
underlying Poisson manifold: it
measures the non invariance of a given volume form with respect to hamiltonian transformations. The vector field (and so the integrated function) depends on the choice of
a volume form but its
cohomology class is independent. We require that the modular function is integrable in the usual dynamical sense
but the hamiltonians in involution must be compatible with the groupoid structure in such a way that the contour level sets inherits the
structure of topological groupoid.
The bihamiltonian system on the projective space provides us with such a system: the modular vector field with respect to the symplectic volume
form is the first hamiltonian vector field of the fundamental Lenard hierarchy.
The hamiltonians can be lifted to the symplectic groupoid and give the multiplicative integrability of the modular function: the procedure is
general but, in the form stated in \cite{BCQT}, it requires that the eigenvalues are global smooth functions. This is true in the projective case but not in
the general Grassmannians. This problem needs a more intrinsic understanding of the polarization and will be addressed in a separate publication.

Let us briefly describe the content of the paper. Let $M_\phi$ be a compact hermitian symmetric space that we see as a $G$-hamiltonian space (let $\g={\rm Lie} G$); $\phi$
denotes the non compact root of the Dynkin diagram of $\g$ associated to the symmetric space.
Our strategy for diagonalizing the
Nijenhuis tensor $N_\phi$ consists first in proving Proposition \ref{reduction_of_momentum_map},
where we show that the eigenvalues of every matrix valued function $\cal M$ solving the {\it master equation} 
$$
N^*_\phi d{\cal M} = d{\cal M}^-{\cal M} + {\cal M} d{\cal M}^+  + r d{\cal M}\;,
$$
define Nijenhuis eigenvalues. See the statement of Proposition \ref{reduction_of_momentum_map} for the explanation of symbols. In Theorem \ref{thm_main} 
we introduce the basic solution of the master equation given by the moment map $\mu$ of the $\g$-action in a representation $R$ that is $\phi$-decomposable
(see Definition \ref{rep_decomposable_root}).
This representation
can be chosen as the fundamental representation in all cases but for $M_\phi=SO(n+2)/SO(n)\times SO(2)$, where we have to choose the spin representation. 
Since $M_\phi$ is a $G$-adjoint orbit, its eigenvalues
are constant and we don't get Njienhuis eigenvalues directly from it. Nevertheless, we get the non trivial solutions to the master equation by a reduction procedure.
Indeed, in Subsection \ref{sec_Rtacos} we introduce case by case a chain of nested subalgebras
\begin{equation}\label{chain_subalgebras_general}
\g \supset\g_1 \supset \g_2\ldots \supset\g_n=0
\end{equation}
together with representation $R_k$ of $\g_k$, such that the moment map of $\g_k$ in the representation $R_k$ solves the master equation.
Theorem \ref{thm_main} relies on an explicit form of the contravariant connection that
encodes the Poisson structure of vector bundles associated to the $G$-principal bundle on $M_\phi$.

In order to show that the obtained eigenvalues are all and that the Nijenhuis operator is of maximal rank the essential ingredient
is the concept of {\it collective complete integrability}. This is a method developed in \cite{GS,GS2,GSbook} for constructing
integrable models. One can consider the algebra of collective hamiltonians $F(\g_1,\ldots\g_n)$ generated by the invariant functions on $\g_i^*$ pulled back through the
moment map. They are in involution and, if the above chain of nested subalgebras satisfies the hypothesis of Proposition \ref{Thimm_method}, define an
integrable model. The most famous integrable model of this form is the Gelfand-Tsetlin model on flag manifolds.
The last step is then to prove integrability of the collective hamiltonians associated to the chain (\ref{chain_subalgebras_general}). 
Since the Nijenhuis eigenvalues are a specific choice of action variables for these integrable model,
the Nijenhuis tensor is of maximal rank.
We call the image of the Nijenhuis eigenvalues the {\it bihamiltonian polytope}. We determine these polytopes case by case.
Moreover let us stress that, thanks to the
bihamiltonian description, the collective hamiltonians are a commutative algebra
also with respect to the Bruhat-Poisson structure.
When $M_\phi$ is the Grassmannian $Gr(k,n)$ then we get the Gelfand-Tsetlin model whose integrability is well established since \cite{GS}. To the best of our knowledge,
in the other cases $M_\phi=Sp(n)/U(n), SO(2n)/U(n), SO(n+2)/SO(n)\times SO(2)$ we get new integrable models and so considerable time is spent in proving integrability and
describing the image of the moment map. This
is the content of Theorems \ref{main_theorem_sun}, \ref{main_theorem_spn}, \ref{main_theorem_so_u} and \ref{main_theorem_so_so}.

The plan of the paper is the following. In Section 2 we recall basic facts about Poisson geometry, Poisson-Lie groups; in particular we recall the notion of Poisson
vector bundle that will be an important tool in our proof. We recall basic notions of Poisson-Nijenhuis structures as well and we briefly sketch the construction
of collective integrable models. In Section 3 we recall basic facts of compact hermitian symmetric spaces and fix notations. In Section 4 we define the Bruhat-Poisson structure. In
Section 5 we give an explicit expression of the contravariant connection defined on associated vector bundles. In Section 6 we introduce the Poisson Nijenhuis structure and develop
the tools needed for the diagonalization. We introduce the master equation in Proposition \ref{reduction_of_momentum_map} and prove that the moment map $\mu$ solves it in Theorem
\ref{thm_main}. Finally, in Subsection \ref{sec_Rtacos} we introduce the chain of subalgebras giving the non trivial solutions of the master equation.
The proof that the collective hamiltonians associated to the chains of subalgebras define a completely integrable model is left to Section \ref{Grassmannians} for $Gr(k,n)$,
Section \ref{sec_sp} for $Sp(n)/U(n)$, Section \ref{sec_so2n} for $SO(2n)/U(n)$ and Section \ref{sec_soodd} for $SO(n+2)/SO(2)\times SO(2)$.

\bigskip
\bigskip
\noindent{\bf Notations}. We will denote with $\g$ the compact form of a complex simple Lie algebra. Let $\t=\t_\C\cap\g$,
where $\t_\C\subset\g_\C$ is a choice of the Cartan subalgebra; let $\Phi$ denote the roots and $\g_\alpha$ with $\alpha\in\Phi$ be the root space.
Let $\Phi^\pm$ be a choice of positive (negative) roots and
$\Pi=\{\alpha_1,\ldots\}$ denote the simple roots.
We denote with $\t^*_+$ the fundamental Weyl chamber.
When we consider the classical cases $\g=\su(n),\so(n),\sp(n)$, we identify $\g$ with
an algebra of matrices and we denote with $f_\g$ the corresponding representation, and we refer to it as the fundamental representation. We denote with
$0_\g$ the one dimensional trivial representation. We denote with $G$ the corresponding matrix group integrating it. 

We recall that a simple root $\alpha_i$ is non compact if the positive roots are all of the form $\alpha = \sum_{j\not=i}c^j\alpha_j$ (of compact type) or
$\alpha=\alpha_i + \sum_{j\not=i} c^j\alpha_j$ (of non compact type). In the following, we list all possible non compact roots.

\begin{figure}[ht]
\begin{picture}(10,50)(-20,-30)
\put(0,0){$A_n$}
\multiput(35,0)(20,0){5}{\circle{8}}
\multiputlist(45,0)(20,0)
{{\line(1,0){12}},{\dashline{2}(-5,0)(5,0)},{\dashline{2}(-5,0)(5,0)},{\line(1,0){12}}}
\multiputlist(35,10)(20,0){$\scriptscriptstyle \alpha_1$,
$\scriptscriptstyle \alpha_2$,$\scriptscriptstyle \alpha_i$,
$\scriptscriptstyle \alpha_{n{-}1}$,$\scriptscriptstyle \alpha_n$}
\put(73,-15){$\scriptstyle \uparrow$}
\end{picture}
\begin{picture}(10,50)(-200,-30)
\put(0,0){$B_n$}
\multiput(35,0)(20,0){4}{\circle{8}}
\put(115.5,0){\circle*{8.5}}
\multiputlist(45,0)(20,0)
{{\line(1,0){12}},{\dashline{2}(-5,0)(5,0)},{\line(1,0){12}}}
\multiputlist(35,10)(20,0){$\scriptscriptstyle \alpha_1$,
$\scriptscriptstyle \alpha_2$,$\scriptscriptstyle \alpha_{n{-}2}$,
$\scriptscriptstyle \alpha_{n{-}1}$,$\scriptscriptstyle \alpha_n$}
\put(99,1.5){\line(1,0){13}}
\put(99,-1.5){\line(1,0){13}}
\put(32.5,-15){$\scriptstyle \uparrow$}
\end{picture}

\begin{picture}(10,50)(-20, -40)
\put(0,0){$C_n$}
\multiput(35,0)(20,0){4}{\circle*{8.5}}
\put(115.5,0){\circle{8}}
\multiputlist(45,0)(20,0)
{{\line(1,0){12}},{\dashline{2}(-5,0)(5,0)},{\line(1,0){12}}}
\multiputlist(35,10)(20,0){$\scriptscriptstyle \alpha_1$,
$\scriptscriptstyle \alpha_2$,$\scriptscriptstyle \alpha_{n{-}2}$,
$\scriptscriptstyle \alpha_{n{-}1}$,$\scriptscriptstyle \alpha_n$}
\put(99,1.5){\line(1,0){12.5}}
\put(98.5,-1.5){\line(1,0){13}}
\put(114,-15){$\scriptstyle \uparrow$}
\end{picture}
\begin{picture}(10,50)(-200,-40)
\put(0,0){$D_n$}
\multiput(35,0)(20,0){4}{\circle{8}}
\put(115.5,15){\circle{8.5}}
\put(115.5,-15){\circle{8.5}}
\multiputlist(45,0)(20,0)
{{\line(1,0){12}},{\dashline{2}(-5,0)(5,0)},{\line(1,0){12}}}
\multiputlist(35,10)(20,0){$\scriptscriptstyle \alpha_1$,
$\scriptscriptstyle \alpha_2$,$\scriptscriptstyle \alpha_{n{-}3}$,
$\scriptscriptstyle \alpha_{n{-}2}$,$\scriptscriptstyle $}
\put(107,25){$\scriptscriptstyle \alpha_{n{-}1}$}
\put(110,-25){$\scriptscriptstyle \alpha_{n}$}
\put(99,1.5){\line(4,3){13}}
\put(98.5,-1.5){\line(4,-3){13.5}}
\put(33.5,-15){$\scriptstyle \uparrow$}
\put(125,13){$\scriptstyle \leftarrow$}
\put(125,-17){$\scriptstyle \leftarrow$}
\end{picture}

\begin{picture}(10,40)(-20,-30)
\put(0,0){$E_6$}
\multiput(35,0)(20,0){5}{\circle{8}}
\multiputlist(45,0)(20,0)
{{\line(1,0){12}},{\line(1,0){12}},{\line(1,0){12}},{\line(1,0){12}}}
\multiputlist(35,10)(20,0){$\scriptscriptstyle \alpha_1$,
$\scriptscriptstyle \alpha_2$,$\hspace{10pt}\scriptscriptstyle \alpha_4 $,
$\scriptscriptstyle \alpha_5$,$\scriptscriptstyle \alpha_6$}
\put(70.5,30){$\scriptscriptstyle \alpha_3$}
\put(75,20){\circle{8}}
\put(74.6,4.1){\line(0,1){12.2}}
\put(33,-15){$\scriptstyle \uparrow$}
\put(113,-15){$\scriptstyle \uparrow$}
\end{picture}
\begin{picture}(10,40)(-200,-30)
\put(0,0){$E_7$}
\multiput(35,0)(20,0){6}{\circle{8}}
\multiputlist(45,0)(20,0)
{{\line(1,0){12}},{\line(1,0){12}},{\line(1,0){12}},{\line(1,0){12}},{\line(1,0){12}}}
\multiputlist(35,10)(20,0){$\scriptscriptstyle \alpha_1$,
$\scriptscriptstyle \alpha_2$,$\hspace{10pt}\scriptscriptstyle \alpha_4 $,
$\scriptscriptstyle \alpha_5$,$\scriptscriptstyle \alpha_6$,$\scriptscriptstyle \alpha_7$}
\put(70.5,30){$\scriptscriptstyle \alpha_3$}
\put(75,20){\circle{8}}
\put(74.6,4.1){\line(0,1){12.2}}
\put(133,-15){$\scriptstyle \uparrow$}
\end{picture}
\caption {Dynkin diagrams with the non compact roots marked.}
\end{figure}
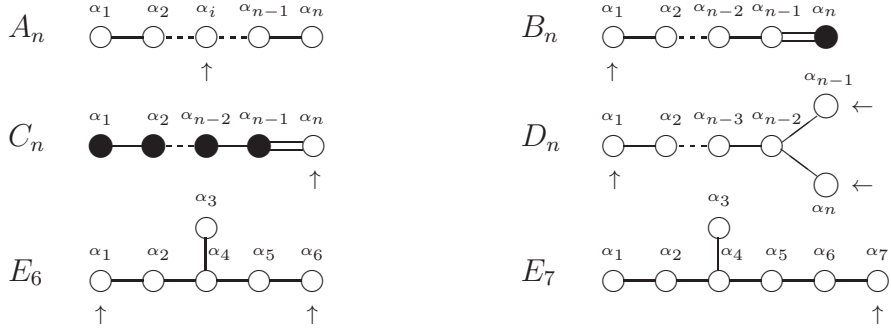

\bigskip
\bigskip
\section{Generalities}

\subsection{Poisson vector bundles}\label{poisson_vector_bundles}
We fix in this Section the conventions and recall basic material about Poisson geometry. Let $(M,\pi)$ be a Poisson manifold with
$\pi\in C^\infty(\Lambda^2TM)$ denoting the Poisson bivector and $\{f,g\}=\pi^{ij}\partial_if \partial_j g$ denoting the Poisson bracket between
$f,g\in C^\infty(M)$. The Jacobi identity for the Poisson bracket can be expressed as $[\pi,\pi]=0$, where $[\ ,\ ]$ denotes the Schouten bracket between
multivector fields. As a consequence the differential $d_{LP}(-)=[\pi,-]$ squares to zero and defines the Lichnerowicz-Poisson cohomology $H_{LP}(M,\pi)$.
Given a volume form $V$ on $M$, the {\it modular vector field} with respect to $V$ is $\chi_V=\divergence_V(\pi)$: it satisfies $d_{LP}(\chi_V)=0$ and
its class in LP cohomology, that does not depend on the choice of the volume form, is called the {\it modular class}.

A Poisson structure defines an algebroid structure on $T^*M$ that we denote as $T^*_\pi M$. The anchor is $\pi:T^*_mM\rightarrow T_mM$, $m\in M$, defined as 
$\langle\pi(\alpha_m),\beta_m \rangle =\langle\pi(m),\alpha_m\wedge\beta_m\rangle$, with $\alpha_m,\beta_m\in T^*_mM$; the bracket on
$\Omega^1(M)$ is defined as
\begin{equation}\label{bracket_algebroid}
\{\alpha,\beta\}_\pi = L_{\pi(\alpha)}\beta - L_{\pi(\beta)}\alpha - d \langle \pi,\alpha\wedge\beta\rangle ~,~~~~\alpha,\beta\in \Omega^1(M) ~~~.
\end{equation}

A Lie group $(G,\pi)$ is called a {\it Poisson-Lie} group if it is a Poisson manifold such that the multiplication is a Poisson map (with
the product Poisson structure on $G\times G$). As a consequence,  $\delta_\g:\g\rightarrow \wedge^2\g$ defined as
$\delta_\g(X)=\frac{d}{dt}\pi(\exp tX)|_{t=0}$,
$X\in\g$, defines a Lie algebra structure on $\g^*$. We call $(\g,\delta_\g)$ a {\it Lie bialgebra}. Let us assume that $G$ is connected and simply connected; let $G^*$
be the connected and simply connected group integrating $\g^*$: it can be shown that
there exists a canonical Poisson-Lie structure on it, such that $(\g^*)^*=\g$ as Lie algebras and $G^*$ is said to be the Poisson-Lie dual of $G$.
The action of  $(G,\pi_G)$ on $(M,\pi_M)$ is a {\it Poisson action} if the action seen as a map from $(G\times M,\pi_G\oplus\pi_M)$ to $(M,\pi_M)$
is a Poisson map. At the infinitesimal level this means
that for each $X\in\g$, denoting with $\ell: X\rightarrow \ell_X$ the map associating the corresponding fundamental vector field on $M$, we have that
$$
L_{\ell_X}(\pi_M) = [\ell_X,\pi_M] = \ell(\delta_{\g}(X)) \;\;.
$$

A subgroup $H\subset G$ is a Poisson-Lie subgroup if it is a Poisson submanifold. Let $\h$ be the Lie algebra of $H$. It can be easily seen that
$H$ is a Poisson-Lie subgroup if and
only $\h^\perp\subset\g^*$ is an ideal of the dual Lie algebra $\g^*$ or equivalently if and only if
$\delta_\g(\h)\subset \wedge^2\h$. Finally, there is a unique Poisson structure on $G/H$ such that the quotient map
$G\rightarrow G/H$ is Poisson.

A vector bundle $\E$ over a Poisson manifold $(M,\pi)$ is a {\it Poisson vector bundle} if there exists a bracket
$\{,\}_\E: C^\infty(M)\otimes \Gamma^\infty(\E)\rightarrow \Gamma^\infty(\E)$, where $\Gamma^\infty(\E)$ denotes the smooth sections, such that, for each
$f,g\in C^\infty(M)$ and $\sigma\in \Gamma^\infty(\E)$ we have
\begin{itemize}
 \item[$i$)] $\{f,g\sigma\}_\E= \{f,g\}\sigma+g\{f,\sigma\}_\E$ ,
 \item[$ii$)] $\{fg,\sigma\}_\E= f\{g,\sigma\}_\E + g \{f,\sigma\}_\E$.
\end{itemize}
See \cite{Ginz1} for a reference. These data can be equivalently encoded in the flat contravariant connection,
$\nabla: \Omega^1(M)\otimes \Gamma^\infty(\E)\rightarrow \Gamma^\infty(\E)$ defined as
$$
\nabla_{df}(\sigma) = \{f,\sigma\}_\E ~, ~~~~~~~ f\in C^\infty(M), \sigma\in\Gamma^\infty(\E) \;.
$$
Another equivalent way of stating the properties of Poisson vector bundle is by saying that $\nabla$ defines a representation of the algebroid $T^*_\pi M$
canonically associated to $(M,\pi)$ (see \cite{ELW} for the definition of an algebroid representation). Let
\begin{equation}\nonumber
 \begin{array}{ccc}
  (P,\pi_P) & \leftarrow & (G,\pi_G) \cr
\downarrow & & \cr
  (M,\pi_M) & &
 \end{array}
\end{equation}
be a {\it Poisson principal bundle}, that is a principal $G$-bundle $P$ over $M$, such that the right action
of $(G,\pi_G)$ on $(P,\pi_P)$ is a Poisson action and the projection $(P,\pi_P)\rightarrow (M,\pi_M)$ is a Poisson map.
Let $R:G\rightarrow {\rm End} V$ be a right representation of $G$ on the vector space $V$ and let $\E_R = V\times_R G$ be the associated vector bundle.
We can characterize sections of $\E_R$ as equivariant functions $C^\infty(P,V)^G$, {\it i.e.} $\sigma\in C^\infty(P,V)^G$ if
$\sigma: P \rightarrow V$ is such that $\sigma(pg)=\sigma(p)R(g)$,
$p\in P$, $g\in G$.

\smallskip
\begin{lemma}\label{poisson_vector_bundle}
The bracket
$$
\{f,\sigma\}_{\E_R} \equiv \{f,\sigma\}_P
$$
between $f\in C^\infty(M)= C^\infty(P)^G$ and $\sigma\in C^\infty(P,V)^G$
endows the associated vector bundle $\E_R$ of the structure of Poisson vector bundle of $(M,\pi_M)$.
\end{lemma}
{\it Proof}.
Since the right $G$ action on $P$ is Poisson, we have that for each $X\in\g$ (denoting with $r: X\rightarrow r_X$
the fundamental vector field of $X\in \g$)
\begin{eqnarray*}
 r_X(\{f,\sigma\}_P) &=& \{r_X(f),\sigma\}_P + \{f,r_X(\sigma)\}_P + \langle r(\delta_{\g}(X)),df\wedge d\sigma\rangle  \cr
&=& \{f,r_X(\sigma)\}_P=\{f,\sigma\}_PR(X) \;,
\end{eqnarray*}
where $\delta_\g$ denotes the bialgebra structure of $\g$; the first term and the third term of the rhs of the first line vanish since $f$ is
invariant with respect to the $\g$ action. \qed

\subsection{Poisson-Nijenhuis structures}\label{Poisson-Nijenhuis structures}

A $(1,1)$ tensor $N:TM\rightarrow TM$ is called a {\it Nijenhuis} tensor if it has vanishing Nijenhuis torsion, {\it i.e.} for any
couple $(v_1,v_2)$ of vector fields on $M$ we have
$$
T(N)(v_1,v_2)=[N v_1,N v_2]-N ([N v_1,v_2] + [v_1,N v_2]-N [v_1,v_2] )=0  ~.
$$
Given any bivector $\pi$, we recall that $\{,\}_\pi$ denotes the antisymmetric bracket on one forms defined in (\ref{bracket_algebroid}).
A triple $(M,\pi,N)$, where $(M,\pi)$ is a Poisson manifold and $N$ a Nijenhuis tensor is called a {\it Poisson-Nijenhuis} (PN) manifold if
$\pi$ and $N$ are compatible, {\it i.e.}
$$
N\circ \pi = \pi \circ N^*~,~~~~ \{\alpha,\beta\}_{N\pi}=\{N^*\alpha,\beta\}_\pi + \{\alpha,N^*\beta\}_\pi - N^*\{\alpha,\beta\}_\pi  ~~~~,
$$
for $\alpha,\beta\in \Omega^1(M)$, where $N^*$ denotes the dual map.

We will consider the case, where $\pi=\Omega^{-1}$ is the inverse of a symplectic form and there exists a compatible
Poisson structure $\pi_0$, {\it i.e.} such that $[\Omega^{-1},\pi_0]=0$, or equivalently there is a pencil of Poisson structures
$\pi_t=\pi_0 + t \Omega^{-1}$ for $t\in\R$. In this case $(M,\Omega^{-1},N=\pi_0\circ\Omega)$ is a Poisson-Nijenhuis structure
(called also $\omega N$-manifold). The PN structures are closely related to integrable systems, see \cite{MM} for a general reference to bihamiltonian systems, here we will recall 
few basic facts.

The {\it spectral problem} associated to the $PN$ structure is the problem of determining the eigenvalues of $N$.
The eigenspace of $N$ corresponding to eigenvalue $\lambda$ is the null space of $\pi_0-\lambda\Omega^{-1}$; since the latter is anti-symmetric, the dimension of the null space
is at least 2. We can then conclude that if $\dim M = 2n$ then $N$ can have at most $n$ distinct eigenvalues.
We say that the rank is maximal if the
distinct eigenvalues are exactly $n$ on a dense open set of $M$. We can define a map $J_N: M\rightarrow \R^n$ associating to every point $m\in M$ the eigenvalues
$(\lambda_1(m),\ldots,\lambda_n(m))\in \R^n$ and we call it the {\it bihamiltonian moment map}. There is not of course a unique way of defining it, according to how we enumerate
the eigenvalues; one possibility is to order them, but in the examples considered in this paper other choices will be more natural.
We call the image of $J_N$ the {\it bihamiltonian polytope} and denote it with
${\cal C}(N)$. Note that for each $t\in\R$, the preimage along $J_N$ of the union of hyperplanes ${\cal C}(t)=\bigcup_k {\cal C}^{(k)}(t)$, where
${\cal C}^{(k)}(t)=\{\lambda\in{\cal C}(N)|\,\lambda_k=-t\}$,
is the set of points where the $\pi_t=\pi_0 +t \Omega^{-1}$ has not maximal rank.

A point $m$ is {\it regular} if ${\rm rk}(dJ_N(m))=n$. If $\{\lambda_i\}$ is a collection of functions that give the eigenvalues of $N$ in
a neighborhood of regular points then they satisfy the following equation
\begin{equation}\label{nijenhuis_eigenvalues}
 N^* d\lambda_i = \lambda_i d\lambda_i ~~~~.
\end{equation}
It can be shown that the eigenvalues $\lambda_i$ are in involution with respect to both $\Omega^{-1}$ and $\pi_0$. A collection of smooth functions
$\{I_k\}_{k\geq 0}$ satisfies the
{\it Lenard recursion relations} if
\begin{equation}\label{lenard_recursion}
d I_{k+1} = N^* d I_k    ~~.
\end{equation}
As a consequence, the $I_k$'s are in involution with respect to both $\Omega^{-1}$ and $\pi_0$. A canonical collection of 
such functions is given by $I_k=\frac{1}{k}\Tr N^k$, $k=1\ldots n$ (this is a consequence of (\ref{nijenhuis_eigenvalues})).

The modular vector field of $\pi_t=\pi_0+t\Omega^{-1}$, $t\in\R$ with respect to the
symplectic volume form is independent on $t$. It is a consequence of
Theorem 3.5 of \cite{DaFe} that this modular vector field is the
$\Omega^{-1}$ hamiltonian vector field of $I_1$, {\it i.e.}
$$
\chi_{\Omega} = \divergence_{\Omega}\pi_t= \Omega^{-1} (d {\rm Tr} N)   ~.
$$
In general, $\chi_\Omega$ is only a Poisson vector field with respect to $\pi_t$. It is easy to show that $\log\det (N+t)$ is a local hamiltonian for $\chi_\Omega$
with respect to $\pi_t$ that is defined on all points such that $-t$ is not an eigenvalue of $N$.

\subsection{Collective complete integrability}\label{collective_complete_integrability}
We recall here a general method for constructing integrable models, called Thimm method in \cite{GS2}, which we refer for details (see also \cite{GSbook}).
Let $M$ be an hamiltonian $K$-space with moment map $\Phi: M\rightarrow \k^*$, where $\k={\rm Lie} K$.
An hamiltonian of the form $\Phi^*(c)$ for $c\in C^\infty(\k^*)$ is called {\it collective}. Any $K$-invariant function $f\in C^\infty(M)^K$ Poisson
commutes with collective hamiltonians.

\begin{defn}\label{definition_multiplicity_free}
An hamiltonian $K$-space $(M,\Phi)$ is multiplicity free if one of the following equivalent conditions is satisfied:
\begin{itemize}
 \item[$i$)] the algebra of $K$-invariant functions $C^\infty(M)^K$ is Poisson commutative;
 \item[$ii$)] for each $\alpha\in\k^*$, denoting with $K_\alpha\subset K$ its stability group with respect to the coadjoint action, the action of
$K_\alpha$ on $\Phi^{-1}(\alpha)$ is transitive;
 \item[$iii$)] for each $\alpha\in\k^*$, denoting with ${\cal O}_\alpha$ the coadjoint orbit through $\alpha$, the action of $K$ on
$\Phi^{-1}({\cal O}_\alpha)$ is transitive.
\end{itemize}
\end{defn}

The equivalence between properties $(i-iii)$ is shown in \cite{GS2}. Let us consider the following chain of subalgebras
$$
\k \equiv \k_1 \supset \k_2 \supset \ldots \k_k\supset\k_{k+1}=\{0\}~~~~,
$$
and let us denote with $K_i\supset K_{i+1}$ the corresponding chain of subgroups. Let $p_i: \k^*\rightarrow \k_i^*$ be the map dual to the 
inclusion $\k_i\subset \k$;
it is easy to see that the invariant functions on $\k_i^*$ pulled back to
$M$ with $(p_i\circ\Phi)^*$ form an abelian subalgebra $F(\k_1,\ldots\k_k)\subset C^\infty(M)$ of the Poisson algebra of functions on $M$.

\smallskip
Let us denote with $p_{ij}:\k_i^*\rightarrow \k_j^*$, $j>i$, the dual map of the inclusion of subalgebras. Every coadjoint orbit ${\cal O}\subset\k_i^*$
is a $K_j$-hamiltonian space with moment map $p_{ij}$. The following result is proven in \cite{GS2}.

\smallskip
\begin{proposition}\label{Thimm_method}
$F(\k_1\ldots\k_k)$ defines a completely integrable model if and only if any coadjoint orbit ${\cal O}\subset \k_i^*$ in the image of
$p_i\circ\Phi: M\rightarrow \k_i^*$  is multiplicity free with respect to the
$K_{i+1}$ coadjoint action for each $i=1,\ldots k$.
\end{proposition}

If the subalgebras $\k_i$ are semisimple then action variables for such an integrable model can be defined as follows.
Let $\beta_i:\k_i^*\rightarrow (\t_{i}^*)_+$ be the map
that sends each point of $\k^*_i$ to the unique intersection of its $K_i$-coadjoint orbit with the positive Weyl chamber.
This is a continuous map that is smooth in the preimage of the interior of the Weyl chamber. Let $\{\xi_i\}$ be a basis of integral 
lattice of $\t_i$:
then the variables $\lambda_i = \langle \xi_i,\beta_i\circ \mu_{\k_i}\rangle$ are action variables.

\smallskip
The most important example of this construction is the so called {\it Gelfand-Tsetlin} integrable model on flag manifolds. We will discuss it 
in the case of Grassmannians in Section \ref{Grassmannians}.

\bigskip
\bigskip

\section{Compact hermitian symmetric spaces}\label{compact_hermitian_spaces}

Let us first fix the geometrical setting of compact hermitian symmetric spaces that we will need later, see \cite{wolf}.

Let $\phi\in\Pi$ be a non compact root and $\Phi_c^+$ and $\Phi_{nc}^+$ the positive roots of compact and non compact type. Let $\h_\phi\subset\g$ be
the Lie subalgebra defined as
$$
\h_\phi = \t \oplus_{\alpha\in\Phi_c^+} (\g_\alpha\oplus \g_{-\alpha}) \cap\g.
$$
and let us denote $H_\phi\subset G$ the closed subgroup integrating it. We denote with $Z(\h_\phi)\subset\h_\phi$ the one
dimensional center. Let $\rho_\phi\in Z(\h_\phi)$ be normalized by $\phi(\rho_\phi)=i$. We denote with $\h_\phi^\perp$ the orthogonal
space to $\h_\phi$ with respect to the Killing form. We have that
$$
\h_\phi^\perp = \oplus_{\alpha\in\Phi^+_{nc}}(\g_\alpha\oplus\g_{-\alpha})\cap\g.
$$

By identifying $\g$ with $\g^*$ thanks to the Killing form, $G/H_\phi$ is identified as the adjoint orbit of $\rho_\phi$. This fixes the KKS
symplectic form $\Omega_{kks}$
in such a way that $G$ acts hamiltonially with moment map given by
\begin{equation}\label{momentum_map}
\mu(g) = g \rho_\phi g^{-1}\;\:      g\in G.
 \end{equation}

The automorphism $\sigma_\phi = \Ad_{K_\phi}$, where $K_\phi=e^{\pi \rho_\phi}$, satisfies $\sigma^2_\phi={\rm id}$ so that $K_\phi^2=e^{2\pi \rho_\phi}\in Z(G)$.
We have accordingly that $\h_\phi$ and $\h_\phi^\perp$ are the eigenspaces of $\sigma_\phi$ corresponding to eigenvalue $1$ and $-1$ respectively.

We have also that $K_\phi^2= e^{2ik_\phi} 1$ so that the fundamental representation decomposes as $V_+\oplus V_-$ corresponding to the eigenvalues $\pm e^{ik_\phi}$ of $K_\phi$. 
Let us denote with $n_\pm= \dim V_{\pm}$.
Let $R_\pm$ denote the representations of $H_\phi$ on $V_\pm$.
Let us consider the homogeneous principal bundle $H_\phi\rightarrow G \rightarrow G/H_\phi$ and let $\E_\pm = G \times_{R_\pm}\!\! V_\pm$ denote the
vector bundles associated to $R_\pm$.

For $g\in G\subset M_N(\C)$ we get the corresponding decomposition
\begin{equation*}
g = \left(\begin{array}{cc}
           g^{++} & g^{+-} \cr
           g^{-+} & g^{--}
          \end{array}
\right)   ~,
\end{equation*}
and let
\begin{eqnarray*}
\sigma_+(g) = \left(\begin{array}{c}
                  g^{++} \cr
                  g^{-+}
                 \end{array}
\right) \in M_{N,n_+}(\C) ~,~~~~~~~~~
\sigma_-(g)= \left(\begin{array}{c}
                  g^{+-} \cr
                  g^{--}
                 \end{array}
\right) \in M_{N,n_-}(\C) ~.\nn
\end{eqnarray*}
We define $\idempotent_{\pm}(g)=\sigma_{\pm}(g)\sigma^\dagger_{\pm}(g)\in M_N(\C)$.

\begin{lemma}\label{lem_3.1}
One has $\idempotent_{\pm}^2=\idempotent_{\pm}$. Moreover $\idempotent_++\idempotent_-= 1 _N$ and
\begin{equation}\label{momentummap_versus_idempotent}
\mu= \sigma_+ R_+(\rho_\phi) \sigma_+^\dagger + \sigma_- R_-(\rho_\phi) \sigma_-^\dagger ~.
\end{equation}
\end{lemma}
{\it Proof}. The first assertions follow from $\sigma^\dagger_\pm\sigma_\pm = 1_{n_\pm}$ that follows from $g^\dagger g=1$, the last one is also clear. \qed
\medskip

The idempotents $\idempotent_\pm$ of Lemma \ref{lem_3.1} define the vector bundles $\E_{\pm}$ as ${\rm Im}\,\idempotent_{\pm}\subset\C^N$. Let us discuss the various cases.

\begin{example}\label{Ex_A}  {\rm  {\bf(AIII)}
Let $G=SU(n)$ and let us choose as non compact root the $k$-th root of the Dynkin diagram. If we choose as Cartan subalgebra
the diagonal matrices we then get $H_\phi= S(U(k)\times U(n-k))$ embedded as block diagonal matrices. The symmetric space is the Grassmannian
$Gr(k,n)$ of $k$-vector space inside $\C^{n}$. We have that
\begin{equation*}\label{An}
 \rho_\phi =\frac{i}{n}\left(\begin{array}{cc}
                    (n-k) 1_{k}& 0\cr
                    0 & - k 1_{n-k}
                  \end{array}
\right),~~
K_\phi = e^{i\pi(1-\frac{k}{n})} \left(\begin{array}{cc}
                                  1_{k} & 0\cr
                                  0& - 1_{n-k}
                                  \end{array}
\right)~.
\end{equation*}
Then clearly we get that $R_+= f_{\u(k)}\times 0_{n-k}$ and $R_-=0_k\times f_{\u(n-k)}$, where $f$ denotes the fundamental representation and
$0$ the trivial one; $\E_+$ is the rank $k$ tautological vector bundle over $Gr(k,n)$ and $\E_-$ is the rank $n-k$ tautological vector bundle over
$Gr(n-k,n)\sim Gr(k,n)$. From (\ref{momentummap_versus_idempotent}) we get
$$
\mu = i\frac{n-k}{n} \idempotent_+ - i \frac{k}{n} \idempotent_- = i \idempotent_+ - i \frac{k}{n} ~~~~.
$$

\begin{rem}{\rm
One can equally write $Gr(k,n)=U(n)/U(k)\times U(n-k)$. Then one can pick
\bea \rho_\phi =\frac{i}{2}\left(\begin{array}{cc}
                     1_{k}& 0\cr
                    0 & - 1_{n-k}
                  \end{array}
\right),~\nn\eea
which shortens some calculations.}
\end{rem}
}
\end{example}
\begin{example}\label{Ex_B}{\bf (BDI odd)} {\rm  Let $G=SO(2n+1)$ and let us consider the first root in the Dynkin diagram, being the unique non compact root. 
Let us choose the Cartan subalgebra as $n$ copies of $\so(2)$, each copy of them embedded in diagonal $2$-dimensional block and having zero in the 
first diagonal entry. Then we have that $H_\phi$ is $SO(2n-1)\times SO(2)$ embedded as block diagonal matrices with $SO(2)$ sitting in the lowest right block. 
The hermitian space is the Grassmannian of oriented real
$2$-dimensional subspaces of $\R^{2n+1}$.
We then have
\begin{equation}\label{rho_so_odd}
\rho_\phi = \left(\begin{array}{cc}
                0_{2n-1} & 0 \cr
                0 & \sigma
               \end{array}
\right),~
K_\phi = \left(\begin{array}{cc}
                1_{2n-1} & 0\cr
                0 & -1_2
               \end{array}
\right)~,
\end{equation}
where
$$\sigma = \left(\begin{array}{cc}
                       0 & 1 \cr
                       -1& 0
                      \end{array}
\right)~. $$
Then clearly $R_+ = f_{\so(2n-1)} \times 0_{\so(2)}$ and $R_- = 0_{\so(2n-1)}\times f_{\so(2)}$. The vector bundle $\E_-$ is the rank $2$
tautological vector bundle. From (\ref{momentummap_versus_idempotent}) we get that
$$
\mu = \sigma_- \sigma \sigma_-^\dagger ~,
$$
so that $\mu^2=\idempotent_-$.
}\end{example}
\begin{example}\label{Ex_C}{\bf (CI)} {\rm Let $G=Sp(n)$ be the compact symplectic group (denoted sometimes as $USp(n)=Sp(2n,\BB{C})\cap U(2n)$). 
In this case the only non compact root is the last one in the Dynkin diagram. The algebra is described as
\begin{equation*}
\label{simplectic_group}
\sp(n) = \left\{\left(\begin{array}{cc}
                      A & B \cr
-B^\dagger& - A^t
                     \end{array}\right)~,~ A=-A^\dagger, B= B^t, A,B\in M_n(\C) \right\} ~~~~.
\end{equation*}
and the stability subgroup is $H_\phi=U(n)$. By choosing the Cartan subalgebra $\t=\R^n$ embedded as $a\in\R^n\rightarrow {{\rm diag}(ia,-ia)}$ we see that
$U(n)$ is embedded as
\begin{equation*}
 H_\phi=U(n) = \left\{\left(\begin{array}{cc}
                             X & 0\cr
                             0 & \bar X
                            \end{array}
\right)~,~ X\in U(n) \right\}~~~.
\end{equation*}
We then have
\begin{equation*}
\rho_\phi =  \left(\begin{array}{cc}
                    \frac{i}{2} 1_n & 0 \cr
                    0 & - \frac{i}{2} 1_n
                   \end{array}
\right)~,~
K_\phi = \left(\begin{array}{cc}
                    i\, 1_n & 0 \cr
                    0 & - i\, 1_n
                   \end{array}
\right) ~.
\end{equation*}
The representations $R_+ = f_{\u(n)}$ and $R_-= \bar{f}_{\u(n)}$ and
\begin{equation*}
\mu = \frac{i}{2} \idempotent_+ - \frac{i}{2}\idempotent_- = i \idempotent_+ - \frac{i}{2} ~~~~~.
\end{equation*}
}
\end{example}
\begin{example}\label{Ex_Dn}{\bf (DIII)} {\rm Let $G=SO(2n)$ and let us consider as non compact root the last root in the Dynkin diagram. The subgroup is then $U(n)$ and the
symmetric space $SO(2n)/U(n)$ is the space of orthogonal complex structures on $\R^{2n}$. Let us choose as Cartan subalgebra
$$
\t = \{\left(\begin{array}{cc}
                0_n& a \cr
                -a & 0_n
               \end{array}
\right)\ , \ a = {\rm diag}(a_1,\ldots, a_n), a_i\in\R\}~.
$$
The subgroup $H_\phi = U(n)$ is then embedded as
$$
A+iB \in U(n) \rightarrow \left(\begin{array}{cc}
                                 A & B \cr
                                 -B & A
                                \end{array}
\right)~.
$$
We then have
\begin{equation*}
\rho_\phi = \left(\begin{array}{cc}
                   0_n & \frac{1}{2} 1_n \cr
                  -\frac{1}{2} 1_n & 0
                  \end{array}
\right)~,~
K_\phi = \left(\begin{array}{cc}
                   0_n &  1_n \cr
                  - 1_n & 0
                  \end{array}
\right)~.
\end{equation*}
The eigenspaces $V_\pm = \langle (a,\pm i a), a \in \C^n \rangle$. By direct computation we see that $R_+ = f_{\u(n)}$ and $R_-= \overline{f}_{\u(n)}$.
By a direct computation we see that
$$
\mu = \frac{i}{2} \idempotent_+ - \frac{i}{2} \idempotent_- = i \idempotent_+ - \frac{i}{2}
$$
}
\end{example}
\begin{example}\label{Ex_D1}{\bf (BDI even)} {\rm Let $G=SO(2n)$ and let us consider the first root in the Dynkin diagram as the non compact root. The subgroup is then
$H_\phi=SO(2(n-1))\times SO(2)$. Let us choose now the Cartan as $\t = \oplus_{k=1}^n \so(2)$ embedded as a diagonal $2\times 2$ block matrix.
We then have
$$
\rho_\phi = \left(\begin{array}{cc}
                0_{2n-2} & 0 \cr
                0 & \sigma
               \end{array}
\right),~
K_\phi = \left(\begin{array}{cc}
                1_{2n-2} & 0\cr
                0 & -1_2
               \end{array}
\right)~.
$$
Analogously to the (BDI odd) case, we have that that $R_+ = f_{\so(2(n-1)}\times {\rm 0}_{\so(2)}$ and $R_-={\rm 0}_{\so(2(n-1))}\times f_{\so(2)}$ and
$\mu = \sigma_- \sigma \sigma_-^\dagger$.
}
\end{example}

\bigskip

\section{The Bruhat-Poisson structure}

We recall here the definition of the Bruhat-Poisson structure on compact hermitian symmetric spaces $G/ H_\phi$. It is obtained from the so called standard
Poisson structure on $G$, that we are going to define first.

Let $G$ be the compact form of the
complex classical group $G_\C\subset SL(N,\C)$, $\g$ and $\g_\C$ be their Lie algebras. Recall that $\g=\{X\in\g_\C ~| ~ X^\dagger=-X\}$.
Let us fix a Cartan subalgebra $\t\subset\g$ and the set of simple roots $\Pi$; we denote with $\Phi^\pm$ the positive (negative)
roots. For each root $\alpha$ we denote with $\g_\alpha\subset\g_\C$ the root space. Let us define $J: \g_\C \rightarrow \g_\C$ as
\begin{equation}
\label{complex_structure_lie_algebra}
J(\t) = 0 ~~,~~~   J(E_\alpha) = \pm i E_\alpha~~~\alpha\in \Phi^\pm~,
\end{equation}
where $E_\alpha\in \g_\alpha$.
Let us remark that if $\h_\phi$ denotes the subalgebra associated to the non compact root $\phi$, as described in the previous section, we have that
\begin{equation}\label{action_rho_perpendicular}
J|_{\h_\phi^\perp} = \ad_{\rho_\phi}~.
\end{equation}

Let us define
\begin{equation}
\label{C+-}
C_{\pm} = i \pm J\;.
\end{equation}
The Iwasawa decomposition is defined as
$\g_\C=\g\oplus \b_\pm$, where $\b_\pm= C_\pm(\g)$; $\g$ and $\b_\pm$ are lagrangian subalgebras with respect to the non degenerate pairing
$\langle A,B\rangle = \Im\Tr[AB]$. The triple $(\g_\C,\g,\b_\pm)$ is a {\it Manin triple}. Let us denote with $(\prg_+,\prgsplus)$ and $(\prg_-,\prgsminus)$ 
the projections defined by the decomposition $\g\oplus\b_+$ and $\g\oplus\b_-$ respectively.
We get in particular an identification of $\g^*$ with $\b_\pm$.
If we use $\Tr$ to identify $\g^*$ with $\g$
then one can check that $C_\pm:\g\rightarrow \b_\pm$ connects these two realizations of $\g^*$.

\begin{example}{\rm
If $\g=\su(n)$ then $\b_\pm= \a\oplus\n_\pm$, where $\a$ denotes the algebra of real diagonal matrices and $\n_+$ ($\n_-$) the strictly upper
(lower) diagonal complex matrices.
The isomorphism
$C_+: \su(n)\rightarrow \b_+$ reads
\begin{equation}\label{isomorphism_dual}
C_+(X)_{rs} =\left\{\begin{array}{ll}
                  2 i X_{rs} & r < s \cr
                  iX_{rr}  & r=s \cr
                  0  & r > s
                 \end{array}
\right.~,~~~~~~~~~~~~~~~ X\in\su(n)
\end{equation}
and analogously for $C_-=C_+^\dagger$. \qed
}
\end{example}

The standard Poisson-Lie structure $\pi_G$ on $G$ (see \cite{LuWe,LuTh} for a general reference) is
defined as
\begin{eqnarray}
\label{standard_poisson}
\langle r_{g^{-1}}\pi_G(g),\xi\wedge\eta\rangle &=& \bra \prg_-(\Ad_{g^{-1}}C_-(\xi)),\,\prgsminus(\Ad_{g^{-1}}C_-(\eta))\ket \cr
                                              &=& -\bra\prg_+(\Ad_{g^{-1}}C_+(\xi)),\,\prgsplus(\Ad_{g^{-1}}C_+(\eta))\ket \;,
\end{eqnarray}
where $C_\pm$ are defined in (\ref{C+-}) and $\xi,\eta\in \g\equiv\g^*$. It can be shown that $\pi_G$ defines a Poisson-Lie structure.
According to the different descriptions of $\g^*$ described above, the dual Lie algebra can be described as the subalgebra $\b_+\equiv\b_-^{op}\subset\g_\C$ or
as the following bracket on $\g$
\begin{equation}
\label{dualliealgebra}
[X,Y]_{\g^*} = [J(X),Y]+[X,J(Y)]  ~,~~~ X,Y\in\g,
\end{equation}
where $J$ is defined in (\ref{complex_structure_lie_algebra}). The dual Poisson-Lie group is the subgroup
$B_+\subset G_\C$ integrating the Lie algebra $\b_+=\b_-^{op}$.
The Iwasawa decomposition
of $G_\C$ consists in the global decomposition $G_\C=G B_+= G B_-$ and
defines the left (right) dressing transformation of $G^*$ on $G$, $(\gamma,g)\rightarrow {}^\gamma g$ and
$(g,\gamma)\rightarrow g^\gamma$, where $g\in G, \gamma\in G^*$, as follows:
\begin{equation}
\label{dressing}
\gamma g={}^\gamma g \gamma^g\,,\;\;\;\;\; g \gamma = {}^g\gamma g^\gamma\;\;\; .
\end{equation}

From the definition of the Poisson bivector $\pi_G$ one can show the following expression for the dressing vector field associated
to $\xi\in\b_+$
\begin{equation}\label{dressing_vector_field}
s_\xi(g) = \frac{d\ }{dt} ({}^{e^{t\xi}}g)|_{t=0} = -\pi_G(r^*_{g^-1} \xi) \;\;\;.
\end{equation}

\smallskip

\begin{lemma}\label{dressing_relations}
The matrix adjoint $\dagger: B_+\rightarrow B_-^{op}$, satisfies  
\begin{equation}
\label{lower_diag_dressing}
 {}^{\gamma^\dagger}g = {}^{\gamma^{-1}}g  \;\;\;,
\end{equation}
for each $g\in G$ and $\gamma\in B_+$. The fundamental vector field of the left dressing action of $\xi\in\b_+$ is
\begin{equation}
\label{fundamental_dressing}
s_\xi(g) = \xi g - g \xi_g = g \xi_g^\dagger - \xi^\dagger g = - s_{\xi^\dagger}(g) \;,
\end{equation}
where $\xi_g=\Ad^*_{g^{-1}}\xi$.
\end{lemma}
{\it Proof}. Since $C_+^\dagger = C_-$, the matrix adjoint sends $\b_+$ in $\b_-$; the statement for the groups follows because they are exponential groups. 
From the above definition of dressing transformation, we get
$${}^{\gamma^\dagger}g (\gamma^\dagger)^g=\gamma^\dagger g = (g^{-1}\gamma)^\dagger=({}^{g^{-1}}\gamma (g^{-1})^\gamma)^\dagger$$
so that ${}^{\gamma^\dagger}g = ((g^{-1})^\gamma)^{-1}=g^{{}^{g^{-1}}\gamma}$, where the last equality follows from $1=(g g^{-1})^\gamma$. Analogously, we apply
the same rules to exchange twice the order of $g{}^{g^{-1}}\gamma$ and find that
$g^{ {}^{g^{-1}}\gamma}= {}^{\gamma^{-1}}g$, from which we get (\ref{lower_diag_dressing}).

We see that
$$\xi g = \frac{d\ }{dt} (e^{t\xi}g)|_{t=0} = \frac{d\ }{dt}( {}^{e^{t\xi}}g)|_{t=0}+ g\frac{d\ }{dt}(e^{t\xi})^g|_{t=0}= s_\xi(g) + g \Ad^*_{g^{-1}}\xi\;,$$
from which the first equality of (\ref{fundamental_dressing}) follows. In the last step we used the fact that the coadjoint action is the derivative of the dressing 
action at the identity. 
The second equality comes by using (\ref{lower_diag_dressing}).
\qed

\medskip
\medskip

Let us consider now a non compact root $\phi$ and let $\h_\phi$ be the subalgebra associated to it and $H_\phi\subset G$ be the subgroup integrating it,
as described in Section \ref{compact_hermitian_spaces}.

\begin{lemma}
The subgroup $H_\phi\subset G$ is a Poisson-Lie subgroup.
\end{lemma}
{\it Proof}. We have to show that $\h^\perp_\phi\subset\g\equiv \g^*$ is an ideal of the Lie bracket (\ref{dualliealgebra}).
Since
$$ [\h_\phi,\h_\phi^\perp]\subset\h_\phi^\perp,~ [\h_\phi^\perp,\h_\phi^\perp]\subset \h_\phi$$
and $J$ preserves $\h_\phi$ and $\h^\perp_\phi$, it is enough to check that $\h^\perp_\phi$ is an (abelian) subalgebra of $\g^*$. Indeed, let
$E_\alpha,E_\beta\in(\h^\perp_\phi)_\C$ be root vectors, then
$$
[E_\alpha,E_\beta]_{\g^*}= i ({\rm sign}\alpha+ {\rm sign}\beta)[E_\alpha,E_\beta] =0\:,
$$
because if $\alpha$ and $\beta$ are both non compact positive (or negative) roots then $\alpha+\beta$ is not a root. \qed

\smallskip

A Poisson structure is then induced on $G/H_\phi$ that we will denote
as $\pi_0$. The quotient map $(G,\pi_G)\rightarrow (G/H_\phi,\pi_0)$ is a Poisson map and the homogeneous $G$ action on $G/H_\phi$ is a Poisson action.
By applying Lemma \ref{poisson_vector_bundle} we conclude that the associated bundles $\E_\pm$ are
Poisson vector bundles, or alternatively that there exists a flat contravariant connection. We will discuss an explicit formula for this connection
in the next section.

\bigskip
\bigskip

\section{The contravariant connection}

Let $M_\phi=G/H_\phi$ denote the compact hermitian symmetric space associated to the non compact simple root $\phi$ (see Section \ref{compact_hermitian_spaces}
for notations). We have seen at the end of the previous section that, if we consider the Bruhat-Poisson structure, the vector bundles
$\E_\pm$ are Poisson vector bundles. In this section we are going to describe
an explicit formula for their contravariant connection.

Let $\nabla: C^\infty(\E_\pm)\rightarrow C^\infty(T M_\phi\otimes\E_\pm)$ be the flat contravariant connection that we define for later convenience as 
$\nabla_{df}(\sigma)= -\{f,\sigma\}$, with $f\in C^\infty(M_\phi)$ and $\sigma:G\rightarrow V_\pm$ equivariant, {\it i.e.} with the opposite sign with respect to Subsection
\ref{poisson_vector_bundles}.
Let us define
$\nabla_{N_\phi}=\Omegakks\circ \nabla: C^\infty(\E_\pm)\rightarrow C^\infty(T^*M_\phi\otimes \E_\pm)\equiv \Omega^1(\E_\pm)$, where $\Omega_{kks}$
is the Kirillov-Konstant-Souriau symplectic form determined by the identification of $M_\phi$ with the adjoint orbit of $\rho_\phi$. The label $N_\phi$
stands for the Nijenhuis tensor to be introduced later in Subsection \ref{def_nijenhuis}.

We recall that $\mu$ is the moment map of the hamiltonian $G$-action defined in (\ref{momentum_map}).
We recall the notations given in Section \ref{compact_hermitian_spaces}. Let $g\in G\subset M_N(\C)$ be written as
$
g=\left(\sigma_{+},\sigma_- 
        \right)
$
with $\sigma_\pm\in M_{N,n_\pm}(\C)$. If we denote with the same symbol the map $g:G\rightarrow M_N(\C)$, we see that the $i$-th row $g_i:G\rightarrow \C^N$ is 
equivariant with respect to the right $H_\phi$
multiplication and defines a section of the trivial vector bundle $\E_+\oplus\E_-$. Analogously,
$(\sigma_\pm)_i$ denotes the $i$-th row and defines a section of $\E_\pm$. 

The main result is given the following proposition.

\smallskip

\begin{proposition}\label{contravariant_connection}
The flat contravariant connection on the trivial bundle $\E_+\oplus \E_-$ reads as
\bea
\nabla_{N_\phi}(g) = ( -J (d\mu) + [\mu,d\mu]) g \ , \;\;\;\; g=\left(\sigma_+,\sigma_-\right)\;.
\label{contr_conn}\eea
Moreover, in the cases (AIII, CI, DIII), the above formula implies for $\E_\pm$
$$
\nabla_{N_\phi}(\sigma_\pm) = \mp C_\pm(d\mu) \sigma_\pm  \;.
$$
\end{proposition}
{\it Proof}. Let us compute $\nabla_{df}(g)=\iota_{df}(\nabla(g)) =-\{f,g\}_{G}$ for $f\in C^\infty(M_\phi)\subset C^\infty(G)$. We see that
$$
\{f,g\}_G = \pi_G(df)(g)= - s_{\xi_f}(g)~~~~~~~
$$
where $\xi_f:G\rightarrow \g^*$ is defined as $\xi_f(g)=r^*_g df$ and the expression of the dressing transformation is given in (\ref{dressing_vector_field}). It is easy to check that $\xi_f(gh)=\xi_f(g)$, for each $h\in H_\phi$ and
$\xi_f'(g)\equiv l^*_g df(g)=\Ad^*_{g^{-1}}\xi_f(g)\in\h^\perp_\phi$. By using formula (\ref{fundamental_dressing}) in Lemma
\ref{dressing_relations}, we see that
\begin{eqnarray*}
s_{\xi_f}(g) &=& -C_-(\xi_f) g + g C_-(\xi_f') = C_+(\xi_f)g - g C_+(\xi_f')\cr
             &=& [-C_-(\xi_f) + g C_-(\xi_f') g^{-1}] g = [C_+(\xi_f) - g C_+(\xi_f')g^{-1}]g ~~~~~~.
\end{eqnarray*}
We have to characterize $\xi_f(g)$ and $\xi_f'(g)$. Since the ring of function of $M_\phi$ is generated by the matrix elements of the moment map $\mu$,
it is enough to consider $f=\mu(X)$, for any $X\in\g$. We are going to show that
$$
\xi_{\mu(X)}= r^*_gd\mu(X) = \{\mu(X),\mu\}_{kks} = - \langle d\mu, v_X \rangle \in \g^*\equiv\g~~~~~~~~,
$$
where $v_X$ is the fundamental vector field of $X$. Indeed, let us evaluate both sides of the above equation with $Y\in\g$: it is easy to check that the result is
$\mu([X,Y])$ on both sides. Analogously, we evaluate
$$
\xi_{\mu(X)}' = [\rho_\phi, g^{-1} X  g] \in \g\equiv \g^*\;\; .
$$
Using (\ref{action_rho_perpendicular}), we then get
\begin{eqnarray*}
g C_\pm(\xi_{\mu(X)}')g^{-1} & = & g ([i\rho_\phi, g^{-1} X  g] \pm [\rho_\phi,[\rho_\phi,g^{-1} X  g]])g^{-1} \cr
                             & = & i [\mu, X] \pm [\mu,[\mu,X]] = \langle - id\mu \mp [\mu,d\mu ],v_X\rangle\; .
\end{eqnarray*}
By collecting all terms and recalling that $v_X=-\Omega_{kks}^{-1}(d\mu(X))$, we get 
$$\nabla_{d\mu(X)} g = \langle \Lambda,v_X\rangle g = - \langle \Lambda,\Omega^{-1}_{kks}(d\mu(X)) \rangle g= \langle d\mu(X),\Omega_{kks}^{-1}(\Lambda)\rangle g\;,$$ 
where
$$\Lambda= \pm ( - C_\pm(d\mu) + i d\mu \pm [\mu,d\mu]) \in \Omega^1 (X_\phi)\otimes\g \;.$$

In the cases (AIII, CI, DIII), observe that $d \mu= \pm i d\idempotent_\pm$ and that $\sigma_\pm^\dagger \sigma_\mp=0$. The result then
follows from an easy computation. \qed

\bigskip
\bigskip

\section{The bihamiltonian system}
\label{diagonalization_nijenhuis}

We recall here the definition of the bihamiltonian system on compact hermitian symmetric spaces and discuss the diagonalization of the Nijenhuis tensor

\subsection{Definition of the Poisson-Nijenhuis structure}\label{def_nijenhuis}
It was proved in \cite{KRR} that the Bruhat and the KKS Poisson structures on compact hermitian symmetric spaces are compatible.
The following argument can be found in \cite{Foth}. The $G$ action on $M_\phi = G/H_\phi$ is Poisson with respect to the Bruhat-Poisson structure
and leaves $\Omega_{kks}^{-1}$ invariant, so that, if we denote with $v:\g\rightarrow \Gamma(TM_\phi)$ the map that associates to $X\in\g$ the fundamental 
vector field $v_X$, we see that
$$
L_{v_X}[\pi_0,\Omega^{-1}_{kks}]= [L_{v_X}\pi_0,\Omega^{-1}_{kks}] = [v(\delta_\g(X)),\Omega^{-1}_{kks}] = 0   ~.
$$
For compact hermitian spaces, the only $\g$-invariant $3$-vector 
field is $0$, so that we conclude that the Poisson structures $\pi_0$ and $\Omega^{-1}_{kks}$ are compatible, {\it i.e.} they satisfy
$$
[\pi_0,\Omega^{-1}_{kks}]=0~~~~~~~~.
$$
The following are direct and fundamental consequences of this fact: $i$) there is a pencil of homogeneous Poisson structures $\pi_t=\pi_0+t\Omega^{-1}_{kks}$, $t\in\R$,
on $M_\phi$. $ii$) The ($1,1$) tensor
\begin{equation}\label{Nijenhuis_tensor}
N_\phi=\pi_0\circ\Omegakks: T M_\phi\rightarrow T M_\phi
\end{equation}
is Nijenhuis, {\it i.e.} it has vanishing Nijenhuis torsion, so that $(M_\phi,\Omega_{kks}^{-1}, N_\phi)$ is a PN structure.

\subsection{Diagonalization of the Nijenhuis tensor}
The main difference between the (AIII, CI, DIII) and the BDI cases is that in the latter case the moment map is not a linear combination of the idempotents
defining the vector bundles $\E_\pm$. This is essentially due to the fact that in the decomposition
of the fundamental representation of $\g=\so(n)$ in eigenspaces of $\exp\pi\rho_\phi$ in the BDI case we
get a reducible representation of $\h_\phi$ where $\rho_\phi$ is not multiple of the identity.
Since this fact plays a central role in our diagonalization of
the Nijenhuis tensor, we have to consider the moment map in a representation where the decomposition is in irreducible components.

Let us consider now an representation $R$ of $\g$ on $V_R$.  
Let $V_R= V_{R^+_\phi}\oplus V_{R^-_\phi}$
be the decomposition in eigenspaces of
$R(e^{\pi\rho_\phi})$; let us call $R^\pm_\phi$ the corresponding representations of $\h_\phi$.
\begin{defn}\label{rep_decomposable_root}
The representation $R$ is decomposable with respect to the non compact root $\phi$ if
\begin{equation}\label{phi_irreducibility}
R^\pm_\phi(\rho_\phi) = r_\phi^\pm 1_{V_{R^\pm_\phi}} ~~~~.
\end{equation}
\end{defn}
\smallskip
It is easy to check that, since $\ad_{\rho_\phi}|_{\h_\phi^\perp}^2=-{\rm id}$, $(r^+_\phi-r^-_\phi)^2=-1$ so that we can choose $V_{R^\pm_\phi}$
such that $r^+_\phi-r^-_\phi=i$.

\begin{example}\label{example_representation}{\rm
We analyze this property case by case using the discussion of the examples of Section \ref{compact_hermitian_spaces}.

\begin{itemize}
 \item[(AIII)] The fundamental representation of $\su(n)$ is decomposable with respect to any non compact root $\alpha_k$; in fact it decomposes into the fundamental 
 representation of $\u(k)$ and $\u(n-k)$ so that $r_{\alpha_k}^+=i(n-k)/n$ and $r_{\alpha_k}^-=-ik/n$. See Example \ref{Ex_A}. 
\item[(CI)] The fundamental representation of $\sp(n)$ is decomposable with respect to the unique non compact root, and the resulting $R_{\phi}^{\pm}$ are the 
fundamental and anti-fundamental representation of $\u(n)$ and
$r_\phi^\pm=\pm i/2$. See Example \ref{Ex_C}.
\item[(DIII)] The fundamental representation of $\so(2n)$ is decomposable with respect to the last root of the Dynkin diagram, where $R_{\phi}^{\pm}$ are the fundamental and 
anti-fundamental representation of $\u(n)$ and $r_\phi^\pm=\pm i/2$. See Example \ref{Ex_Dn}.
\item[(BDI)] The fundamental representation of $\so(n)$ is not decomposable with respect to the first root, as can be seen in Examples \ref{Ex_B} and \ref{Ex_D1}. Their 
spin representations are decomposable with respect to the first root of their Dynkin diagram: indeed the weights are $(\pm1/2,\cdots,\pm1/2)$ so that $r_{\alpha_1}^\pm=\pm i/2$. 
We will give additional details at the end of this section.
\end{itemize}

}\end{example}

Let $R$ be $\phi$-decomposable and let $\E_{R^\pm_\phi}=G\times_{R^\pm_\phi}V_{R^\pm_\phi}$ be the vector bundles on $M_\phi$ associated to $R^\pm_\phi$.
By applying to
(\ref{momentummap_versus_idempotent}) the representation $R$ of the simply connected group integrating $\g$  and denoting $\mu_ R\equiv R(\mu)$ the moment map in this representation, we get
\begin{equation}
\label{momentummap_versus_idempotent_special_rep}
\mu_R = r^+_\phi \idempotent_{R^+_\phi} + r^-_\phi \idempotent_{R^-_\phi} = i \idempotent_{R^+_\phi} + r^-_\phi~~~,
\end{equation}
where $\idempotent_{R^\pm_\phi}$ are idempotents defining $\E_{R^\pm_\phi}$, see the definition given in Lemma \ref{lem_3.1}. 

\begin{theorem}\label{thm_main}
Let $R$ be a representation of $\g$ decomposable with respect to $\phi$. We have that
\begin{equation}\label{fundamental_formula}
N_\phi^* d\mu_R = \pm\, d\mu^\mp_R\, \mu_R \pm \mu_R\,d\mu^\pm_R \mp 2i r^\pm_\phi d\mu_R ~,
\end{equation}
where $d\mu^\pm_R\equiv R(C_\pm(d\mu))$.
\end{theorem}
{\it Proof}.  By using the formula (\ref{contr_conn}) for the contravariant connection we see that
$$N_\phi^*d\mu = \nabla_{N_\phi}(g)\rho_\phi g^\dagger + g\rho_\phi \nabla_{N_\phi}(g^\dagger)= [-J(d\mu),\mu] + [[\mu,d\mu],\mu]\;.$$
Let us show that $[[\mu,d\mu],\mu]=d\mu$. Indeed, for each $X\in\g$ we see that
\begin{eqnarray*}
\langle [[\mu,d\mu],\mu],v_X \rangle&=&[[\mu,[\mu,X]],\mu] = -g \ad_{\rho_\phi}^2[\rho_\phi, g^{-1}Xg] g^{-1}=g[\rho_\phi,g^{-1}Xg]g\cr
&=& [\mu, X] = \langle d\mu, v_X\rangle \;.
\end{eqnarray*}
We then see that
\begin{equation}\label{connection_momentum_map}
 N_\phi^*d\mu = [-J(d\mu),\mu] + d\mu~~~.
\end{equation}
By using the definition of $C_\pm=i\pm J$, we write $J(d\mu)$ as $C_+(d\mu)-i\,d\mu$, and place (\ref{connection_momentum_map}) in the representation $R$
\bea N_\phi^*d\mu_R&=&-d\mu^+_R\,\mu_R+\mu_R\, d\mu^+_R+i\,d\mu_R^2-2i\mu_R\, d\mu_R + d\mu_R~\nn\\
&=&-d\mu^+_R\,\mu_R-\mu_R\, d\mu^-_R+i\,d\mu_R^2+ d\mu_R~.\nn\eea
By using (\ref{momentummap_versus_idempotent_special_rep}) and the idempotency of $\idempotent_{R_{\phi}^+}$ we get 
$d\mu_R^2= 2r_{\phi}^-d\mu_R+i\,d\mu_R$, and so the result.

\qed

\smallskip
The equation (\ref{fundamental_formula}) satisfied by $\mu_R$ is the basic tool for producing eigenvalues of $N_\phi^*$ as we are going to show in the 
following proposition. 

\begin{proposition}\label{reduction_of_momentum_map}
Let ${\cal M}$ be a matrix valued function on an open subset $U\subset M_\phi$. Assume also that the eigenvalue $m$ of ${\cal M}$ is a smooth non constant
function on $U$ with constant multiplicity. Consider an equation of type 
\bea N_\phi^*d{\cal M}=d{\cal M}^-\,{\cal M}+{\cal M}\,d{\cal M}^++rd{\cal M},\label{master_equation}\eea
with $d{\cal M}^++d{\cal M}^-=kd{\cal M}$ and $k,r\in\C$.
Then
\bea N_\phi^*dm=\big(m k+r\big)dm,\label{eigenvalues_moment_map_versus_nijenhuis}\eea
i.e. $m k+r$ is an eigenvalue of $N_\phi^*$.
\end{proposition}
{\it Proof.} Let $x\in C^\infty(U)$ and let $P(x,{\cal M})=\det(Ix-{\cal M})$. We have that
\bea N_\phi^*dP=(N_\phi^*dx)\partial_xP-P\,\Tr[(Ix-{\cal M})^{-1}N_\phi^*d{\cal M}]\,,\nn\eea
where we used the formula $d \det A=\det A\,\Tr[A^{-1}d A]$. We use (\ref{master_equation}) and we get
\bea N_\phi^*dP=(N_\phi^*dx) \partial_xP-P\,\Tr\big[(Ix-{\cal M})^{-1}\big(d{\cal M}^-\,{\cal M}+{\cal M}\,d{\cal M}^++rd{\cal M}\big)\big].\nn\eea
We write the first term of the trace as 
\bea \Tr\big[{\cal M}(Ix-{\cal M})^{-1} d{\cal M}^-\big]&=&\Tr\big[({\cal M}-Ix+Ix)(Ix-{\cal M})^{-1}d{\cal M}^-\big]\nn\\
&=&-\Tr\big[d{\cal M}^-\big]+x\Tr\big[(Ix-{\cal M})^{-1}d{\cal M}^-\big].\nn\eea
We do the same for the second term, and the two combine into
\bea N_\phi^*dP=(N_\phi^*dx)\partial_xP+(xk+r)d_{\cal M} P+kP\,\Tr[d{\cal M}]\nn\ ,\eea
where we denote with $d_{\cal M}P$ the differential of $P$ keeping $x$ fixed.
Let $\alpha$ be the multiplicity of $m$ so that $P(x,{\cal M})=(x-m)^\alpha P_0(x,{\cal M})$ with $P_0(m,{\cal M})\not=0$.
Now we evaluate the above equation at $x=m+\varepsilon$.
It is easy to see that the dominant term of order $\varepsilon^{\alpha-1}$ in $\varepsilon\rightarrow 0$ gives the formula (\ref{eigenvalues_moment_map_versus_nijenhuis}).
\qed

\smallskip
We call (\ref{master_equation}) the {\it master equation}. By Theorem \ref{thm_main} $\mu_R$ satisfies the master equation with 
$d{\cal M}^\pm=d\mu^\pm_R$ and $k=2i,\ r=-2ir^+_\phi$. Of course the eigenvalues of $\mu_R$ are constant and Proposition \ref{reduction_of_momentum_map} 
does not apply.
We will see in the following subsection the general strategy to produce the Nienhuijs eigenvalues.

\subsection{Reduction to a chain of subalgebras}\label{sec_Rtacos}

In order to build the eigenvalues of the Nijenhuis tensor $N_\phi$ we will pick a chain of subalgebras 
$$
\g\supset \g_1 \supset \g_2 \supset \ldots \supset\g_k\ ,
$$
where each $\g_i$ is equipped with a representation $R_i$ such that the moment map in this representation $\mu_{\g_i R_i}$ solves the master equation 
(\ref{master_equation}). With these data, we will get the eigenvalues by applying Proposition \ref{reduction_of_momentum_map} at each step.

In this subsection we will show how to define these data case by case. 
The proof that we get all the eigenvalues from this construction is postponed to the next sections where we will use the 
results about integrability of the collective hamiltonians defined by the above chain of subalgebras.

\smallskip

\noindent\emph{\bf AIII.} Let $M_\phi= Gr(k,n)$; from the discussion in Example \ref{example_representation}, we can conclude that 
Equation (\ref{fundamental_formula}) is valid with $R$ being the fundamental 
representation so that $\mu_R=\mu$ and
$$
N^*_{\alpha_k} d\mu = C_-(d\mu) \mu + d\mu C_+(d\mu) -2i r^+_{\alpha_k} d\mu\;,
$$
where $r^+_{\alpha_k}=i(n-k)/n$.
Since $C_+(d\mu)$ and $C_-(d\mu)$ can be chosen as upper and lower triangular matrices respectively, it is easy to check 
that every $(n-s)\times (n-s)$ upper left minor $\mu^{(s)}$ solves the master equation (\ref{master_equation}) with 
$d{\cal M}^\pm=dC_\pm(\mu)^{(s)}$, $k=2i$ and $r=-2ir_{\alpha_k}^+$. 

In order to read these minors as moment maps of a chain of subalgebras, it is better to look at $Gr(k,n)$ as a $\u(n)$ hamiltonian space rather than
$\su(n)$ and consider the chain of subalgebras 
\bea
\u(n)\supset \u(n-1)\ldots \supset \u(1)\label{chain_subalgebra_su}
\eea
with $\g_s= \u(n-s)$ embedded as the upper-left corner of $\g_{s-1}=\u(n-s+1)$. It is clear that the minor $\mu^{(s)}$ is the moment map of
$\u(n-s)$ in the fundamental representation.
 
The eigenvalues of $\mu^{(s)}$
are the classical Gelfand-Tsetlin variables. In Section \ref{Grassmannians} 
we will review their properties and show that they exhaust all the possible Nijenhuis eigenvalues.

\smallskip

\noindent\emph{{\bf CI} and {\bf DIII.}} From the discussion in Example \ref{example_representation} we know that in both cases the fundamental 
representation of $\g=\so(2n),\sp(n)$ is decomposable with respect to the non compact root $\phi$. Equation (\ref{fundamental_formula}) is then 
valid in the fundamental representation with $r^\pm_\phi=\pm i/2$ in both cases.

We pick the chain of subalgebras as 
\bea \sp(n)\supset \u(n)\supset \u(n-1)\cdots\supset \u(1),\label{GC_spn}\\
\so(2n)\supset \u(n)\supset \u(n-1)\cdots\supset \u(1),\label{GC_so2n}\eea
where $\g_k=\u(n+1-k)$ is embedded as upper left block of $\g_{k-1}$ and is considered in the fundamental representation.

We will show first that the master equation is valid for the first $\u(n)$ step. We need the following general discussion. Let $R$ be a representation of $\g$ 
decomposable with respect to the non compact root $\phi$ and let $V_R=V_{R^+_\phi}\oplus V_{R^-_\phi}$. Since
$R(\h^\perp_\phi):V_{R^\pm_\phi}\rightarrow V_{R^\mp_\phi}$, the moment map $\mu=\mu_{\h_\phi} + \mu_{\h_\phi^\perp}$ in the 
representation $R$ accordingly decomposes as 
$$ \mu_R=\left(\begin{array}{cc}
                \mu_{\h_\phi R^+_\phi} & \mu_{\h_\phi^\perp}^{+-} \cr
                \mu_{\h^\perp_\phi}^{-+} & \mu_{\h_\phi R^-_\phi}
               \end{array}
\right)\;,$$
where $\mu_{\h_\phi R^\pm_\phi}$ is the moment map of $\h_\phi$ in the representation $R^\pm_\phi$.
\smallskip
\begin{lemma}\label{representation} Let $R$ be a representation decomposable with respect to the non compact root $\phi$. For each $\phi$-non compact
positive root $\alpha$ we have that
$$R(E_\alpha) V_{R^+_\phi} = R(E_{-\alpha}) V_{R^-_\phi} = 0~~~~.$$
\end{lemma}
{\it Proof}. Let $\alpha$ be a positive $\phi$-non
compact root. We see that for each $v_+\in V_{R^+_\phi}$
\begin{eqnarray*}
r^-_\phi R(E_\alpha)v_+&=& R(\rho_\phi)R(E_\alpha)v_+ =  R(E_\alpha)R(\rho_\phi)v_+ + R([\rho_\phi,E_\alpha]) v_+ \cr
&=& (r^+_\phi + i) R(E_\alpha) v_+= (r^-_\phi + 2i) R(E_\alpha) v_+~,
\end{eqnarray*}
so that $R(E_\alpha)v_+=0$. Analogously one can show that $R(E_{-\alpha})v_-=0$. \qed

\smallskip
As a consequence the matrices representing $C_+(\h^{\perp}_{\phi})$ and $C_-(\h^{\perp}_{\phi})$ are concentrated in the $(+-)$ and $(-+)$ block respectively. This has the consequence that the $(++)$ and $(--)$ block of the equation (\ref{fundamental_formula}) satisfied by $\mu_R$ 
is the master equation for $\mu_{\h_\phi R^\pm_{\phi}}$. If we consider the $(++)$ component in our case of $\g=\so(2n),\sp(n)$, we get that 
the $\u(n)$ moment map $\mu_{\u(n)}$ in the fundamental representation satisfies the master equation
\bea
N_\phi^* d\mu_{\u(n)} = C_-(d\mu_{\u(n)}) \mu_{\u(n)} + \mu_{\u(n)}  C_+(d\mu_{\u(n)}) + d\mu_{\u(n)}\label{temp_IV} \:.
\eea

The subsequent reductions will proceed exactly the same as in the AIII case. In Sections \ref{sec_sp} and 
\ref{sec_so2n}, we will 
carry out the remaining details, including establishing the independence and the range of the eigenvalues.

\smallskip

\noindent\emph{\bf BDI}. This is the case where we have to use the moment map in a representation different from the fundamental. 
As it was observed in Example \ref{example_representation}, the spin representation $S$ is decomposable with respect to the non compact root 
$\phi=\alpha_1$ with 
$r^\pm_{\phi}=\pm i/2$ so that equation (\ref{fundamental_formula}) for $\mu_S$ means 
$$
N^*d\mu_S = d\mu_S^- \mu_S + \mu_S d\mu_S^+ + d\mu_S\;,
$$
where $d\mu_S^\pm= S(C_\pm(d\mu))$. 

Let $\g=\so(n+2)$ where $n+2=2N,2N+1$. Let us recall a few basic facts of the spin representation $S$.
We label coordinates of $\R^{2N}$ as $\{x_i, i=1,\ldots 2N\}$, and that of $\R^{2N+1}$ as
$\{x_0,x_i, i=1,\ldots 2N\}$. We introduce complex coordinates
$z_i=(x_{2i-1}+i x_{2i})/2$, $i=1,\ldots,N$ and gamma matrices $\Gamma_i$. The action of the gamma matrices on
$V_S=\wedge \langle d\bar{z}\rangle_{i=1}^N$ is defined as $\Gamma_i = d\bar{z}_i\wedge$,
$\Gamma_{\bar{i}}=\iota_{\partial_{\bar{z}_i}}$ and $\Gamma_0=(-1)^{\deg}$. Recalling that $S(X)=\frac{1}{8} X_{ij}[\Gamma_i,\Gamma_j]$ for $X\in\g$ 
we easily see
that
\bea S(\rho_\phi) = i(\Gamma_{\bar{N}}\Gamma_N -\frac{1}{2})~,\label{last_rot}
\eea
so that $S(\rho_\phi) |_{V^\pm_S}= r^\pm_\phi 1_{V^\pm_S}$, where $V_S^+ = \wedge\langle d\bar{z}_i. i=1,\ldots,N-1\rangle$ and
$V_S^-=V_S^+\otimes d\bar{z}_N$ and $S^\pm$ is the representation $(S,\pm i/2)$ of $\so(n)\oplus\so(2)$. Pay attention that
$S^\pm$ is not to be confused with the chirality in the even case. 
With our choice of the Cartan subalgebra, the positive root vectors are represented as
$$
\Delta^+_{2N+1}=\{-dz^{\bar j}\,\iota_{d\bar{z}_i},~i>j;~~\iota_{d\bar{z}_i}(-1)^{\deg};~~\iota_{d\bar{z}_i}\iota_{d\bar{z}_j}\},\;
$$
$$
\Delta^+_{2N}=\{-dz^{\bar j}\,\iota_{d\bar{z}_i},~i>j;~~\iota_{d\bar{z}_i}\iota_{d\bar{z}_j}\} \;.
$$
A basis of $V_S$ is given by the words $d\bar z_{i_1}\wedge\cdots\wedge d\bar z_{i_p}\in V_S$, $i_1<\cdots <i_p$. 
We pick an ordering of the words such that $d\bar z_{i_1}\wedge\cdots\wedge d\bar z_{i_p} \prec d\bar z_{j_1}\wedge\cdots\wedge d\bar z_{j_q}$, if $i_p<j_q$, 
or in case $i_p=j_q$ then $i_{p-1}< j_{q-1}$ and so on. 
In this basis, positive root vectors are upper diagonal so that $C_+(\g)$ are 
upper triangular matrices. We can again use the same logic as for the AIII case and 
conclude that every upper left minor of $\mu_S$ satisfies the master equation with $k=2i$ and $r=1$. In particular the upper left $2^{N-1}$ minor 
$\mu_S^{(N-1)}$ is the moment map for the subalgebra $\so(n)\oplus\so(2)$ in the representation $(S,i/2)$. By iterating the procedure we can
conclude that the upper left $2^{N-s}$ minor $\mu_S^{(N-s)}$ is the moment map of
$$\g_s=\so(n+2-2s)\oplus\underbrace{\so(2)\oplus\ldots \so(2)}_s$$
in the representation $(S,i/2,\ldots,i/2)$.

To summarize, denoting with $\t=\so(2)\oplus\so(2)\ldots$ the Cartan subalgebra of $\so(n+2)$, we proved that we produce Nijenhuis eigenvalues 
considering the eigenvalues of the moment map of the subalgebras appearing in the following chain
\bea
\so(n+2) \supset \so(n)\oplus\so(2) \supset \so(n-2)\oplus\so(2)\oplus\so(2)\supset \ldots \supset \t \label{chain_subgroups_so}
\eea
considered in the representation $(S,i/2,\ldots,i/2)$. The proof of their independence and description of their range will 
be given in Section \ref{sec_soodd}.

\section{$M_\phi= Gr(k,n)=SU(n)/S(U(k)\times U(n-k)))$}\label{Grassmannians}
Let us consider $M_\phi= Gr(k,n)=SU(n)/S(U(k)\times U(n-k))$.
We showed in Subsection \ref{sec_Rtacos} that the moment map $\mu_{\u(n-s)}$ of the subalgebra $\u(n-s)$ appearing
in the chain (\ref{chain_subalgebra_su}) in the fundamental representation solves the master equation (\ref{master_equation}). By applying
Proposition \ref{reduction_of_momentum_map} we get the Nijenhuis eigenvalues. The eigenvalues of these moment maps are the so called
Gelfand-Tsetlin variables. Their integrability has been established in \cite{GS2,GS}; let us briefly recall the construction.

The result is a consequence of the following proposition proved in \cite{GS2}.

\smallskip
\begin{proposition}\label{multiplicity_freeness_un}
Let ${\cal O}$ be a coadjoint orbit of $\u(n)$ and let us consider $\u(n-1)\subset\u(n)$ (embedded in the upper left corner, for instance).
Then ${\cal O}$ is multiplicity free as hamiltonian $U(n-1)$-space.
\end{proposition}

By applying Proposition \ref{Thimm_method} we conclude that the chain
\begin{equation}
\label{GC-chain}
\u(n-1)\supset \u(n-2)\ldots \supset \u(1)\supset 0
\end{equation}
defines an integrable model on any $U(n)$ coadjoint orbit. 

Let us consider the $U(n)$ orbit ${\cal O}_{\tilde\lambda}$ of $i\tilde\lambda$, where $\tilde\lambda=(\tilde\lambda_1, \tilde\lambda_2, \cdots , \tilde\lambda_{n})\in \R^{n}$ satisfies 
$\tilde\lambda_1\leq \tilde\lambda_2\leq \cdots\leq  \tilde\lambda_{n}$ and let again $\mu$ be the $\u(n)$ moment map. 

We recall that the moment map $\mu_{\u(n-s)}$ is the upper left $(n-s)\times(n-s)$ minor $\mu^{(s)}$ of $\mu$. It follows from the mini-max principle (see \cite{GS}) that
the eigenvalues $i\tilde\lambda^{(s)}_j$ of $\mu_{\u(n-s)}$ satisfy the Gelfand-Tsetlin inequalities
\begin{equation}\label{GC_inequalities}
 \tilde\lambda_i^{(s)}\leq \tilde\lambda_i^{(s+1)} \leq \tilde\lambda_{i+1}^{(s)}\,,
\end{equation}
with $i=1,\ldots , n-s$ and $\tilde\lambda_i^{(0)}=\tilde\lambda_i$. The {\it Gelfand-Tsetlin polytope} is defined as the subset
${\cal C}_{GC}(\lambda)\subset \R^{N(\tilde\lambda)}$, with $N(\tilde\lambda)=\dim{\cal O}_{\tilde\lambda}/2$,
of independent solutions of the inequalities (\ref{GC_inequalities}). The $\tilde\lambda_i^{(s)}$ are a choice of action variables of the integrable system 
defined by the chain (\ref{GC-chain}).

Here we are interested to the case of the Grassmannian $Gr(k,n)$ where 
$$\tilde\lambda=(\underbrace{-k/n,\ldots}_{n-k},\underbrace{1-k/n,\ldots}_{k})\;,$$
{\it i.e.} the ordered eigenvalues of $\rho_{\phi}$ defined in Example \ref{Ex_A}.

Then $-i\mu^{(1)}$ has only one non-constant eigenvalue $\tilde\lambda_{n-k}^{(1)}\in[-k/n,1-k/n]$. This procedure can be iterated to the subsequent subalgebras, e.g. 
$-i\mu^{(2)}$ has two non-constant eigenvalues $\tilde\lambda_{n-k-1}^{(2)},\tilde\lambda_{n-k}^{(2)}$ within the range 
$-k/n\leq \tilde \lambda^{(2)}_{n-k-1} \leq \tilde\lambda_{n-k}^{(1)}\leq \tilde\lambda_{n-k}^{(2)}\leq 1-k/n$, 
and so on. 

As an example, for $Gr(2,4)$, we have the pattern
\bea
       \begin{array}{ccccccc}
         -\frac12 &  & -\frac12 &  & \frac12 &  & \frac12 \\
                  & -\frac12 & & \tilde\lambda_2^{(1)} & & \frac12 & \\
                  &  & \tilde\lambda^{(2)}_1 & & \tilde\lambda^{(2)}_2 & \\
                  &  &  & \tilde\lambda_1^{(3)} \\       \end{array}\;.\nn
\eea

From \cite{GS} we know that the Gelfand-Tsetlin variables are
independent and define a completely integrable system. As a consequence they exhaust all the possible eigenvalues of the Nienhuijs tensor $N_\phi$.
We have then shown the following result.

\begin{theorem}\label{main_theorem_sun}
The Nijenhuis tensor (\ref{Nijenhuis_tensor}) on $Gr(k,n)$ is of maximal rank and its eigenvalues are written in 
terms of the Gelfand-Tsetlin variables as
\begin{equation}\label{eigenvalues_grassmannians}
\lambda^{(s)}_j = -\, 2 (\tilde{\lambda}^{(s)}_j - \frac{n-k}{n}).
\end{equation}
The bihamiltonian polytope coincides with the Gelfand-Tsetlin polytope
$${\cal C}_{GC}(\underbrace{0,\ldots,0}_k,\underbrace{2,\ldots,2}_{n-k})~~~~.$$
\end{theorem}

\medskip
\begin{rem}
{\rm The case $k=1$, the complex projective plane $\C P_n$, was solved in \cite{Foth}. There is one non constant eigenvalue $\tilde{\lambda}^{(s)}_{n-s}$ for each $\u(n-s)$. 
As it was observed in \cite{BCQT}, these eigenvalues correspond to a specific basis of the Cartan subalgebra. In fact, since $\tilde{\lambda}^{(s)}_{n-s}$ is the unique non 
constant eigenvalue of $\mu^{(s)}$, we have that $id\tilde\lambda^{(n-s)}=d\Tr\mu^{(s)}$. One then checks that
$d\lambda^{(s)}_{n-s}=d\mu(H_s)$ where 
$$H_s=2i\ {\rm diag}(\underbrace{1,\ldots,1}_{n-s},\underbrace{0,\ldots,0}_s)\:.$$
In particular these eigenvalues are global smooth functions. The result for the general case $Gr(k,n)$ was only conjectured in \cite{Foth}. The eigenvalues are only continuous 
functions; by repeating the above logic one can show that for each $s$, $\sum_j\lambda^{(s)}_j$ is $\mu(H_s)$ up to a constant and is in particular smooth. 
}
\end{rem}

\section{$M_\phi= Sp(n)/U(n)$}\label{sec_sp}
Consider now  $M_\phi=Sp(n)/U(n)$. We showed in Section
\ref{sec_Rtacos} that the moment map in the fundamental representation of $\g_k=\u(n+1-k)\subset \sp(n)$ appearing in the chain 
(\ref{GC_spn}) solves the master equation (\ref{master_equation}). By applying Proposition \ref{reduction_of_momentum_map} we define
the Njienhuis eigenvalues. In the following theorem we show that they are independent by proving the complete integrability of the
collective hamiltonians defined by the chain (\ref{GC_spn}).

\begin{theorem}\label{main_theorem_spn}
The collective hamiltonians $F(\u(n)\ldots\u(1))$ define a completely integrable model.
The Nijenhuis tensor (\ref{Nijenhuis_tensor}) on $Sp(n)/U(n)$ has maximal rank. Its eigenvalues are all obtained as
$$
\lambda_i^{(k)} = -\,2 \tilde\lambda^{(k)}_i  + 1,~ i=1,\ldots n+1-k,
$$
where $i\tilde\lambda^{(k)}_i$ are the eigenvalues of the moment map of the hamiltonian $\g_k=\u(n+1-k)$ action.

The image of the bihamiltonian moment map is then described as the following polytope ${\cal C}(N_\phi)\subset\R^{\frac{n(n+1)}{2}}$, where
$\big(\lambda_i^{(k)}\big)_{1\leq i\leq n+1-k\leq n}\in{\cal C}(N_\phi)$ if
$$0\leq\lambda_1^{(1)}\leq\cdots\leq \lambda_n^{(1)} \leq 2\,,\qquad \lambda_i^{(k)}\leq \lambda_i^{(k+1)} \leq \lambda_{i+1}^{(k)}\ .$$
\end{theorem}

{\it Proof}. We want to apply Proposition \ref{Thimm_method}. It is enough to show the multiplicity freeness of $U(n)$-orbits in $M_\phi$; the
other steps involve orbits of $U(k)$ contained in $\mu_{\u(k)}(M_\phi)$ with respect to the $U(k-1)$ action, that are always multiplicity free
as a consequence of Proposition \ref{multiplicity_freeness_un}.
In order to use condition $ii)$ of Definition \ref{definition_multiplicity_free} we shall show that for almost all $i\tilde\lambda\in\t_n$, the
diagonal $n \times n$ matrices of $\u(n)$, the action of $U(n)_{i\tilde\lambda}$ on $\mu_{\u(n)}^{-1}(i\tilde\lambda)$, where
$U(n)_{i\tilde\lambda}\subset U(n)$ is the stability subgroup of $i\tilde\lambda$, is transitive.

If we parametrize $g\in Sp(n)$ as
\begin{equation}\label{parametrization_spn}
g= \left(\begin{array}{cc}
          A & B \cr
          -\bar{B} & \bar{A}
         \end{array}
\right),~~~~   \begin{array}{cc}A^\dagger A + B^t \bar{B} = 1,& A^\dagger B = B^t \bar{A},\cr AA^\dagger + BB^\dagger=1, & AB^t=BA^t,\end{array}
\end{equation}
we compute
$$
X= g \rho_\phi g^{-1} = \left(\begin{array}{cc}
                               i(A A^\dagger - 1/2) & -i A B^t \cr
                               -i \bar{B} A^\dagger & -i (\bar{A} A^t - 1/2)
                              \end{array}
\right)~,
$$
so that the moment map for the $\u(n)$ action is
$$
\mu_{\u(n)}(X) = i(A A^\dagger - 1/2)~,
$$
and
\begin{equation}
\label{temp_VIII}
\mu_{\u(n)}^{-1}(i\tilde\lambda) =\{\Omega = -i A B^t, A A^\dagger = 1/2 +\tilde\lambda,\, A,B \, {\rm satisfying}\, (\ref{parametrization_spn}) \}\:.
\end{equation}

This constraints $1/2\pm \tilde\lambda \geq 0$. It is easy to see that the action of $k\in U(n)_{i\tilde\lambda}$ on $\Omega$ reads
$k \Omega k^t$. Let us define
$$
A_0 = \sqrt{1/2 + \tilde\lambda},\, B_0 = \sqrt{1/2 - \tilde\lambda},
$$
so that $\Omega_0 = -i A_0 B_0 =-i \sqrt{1/4 - \tilde\lambda^2}\in \mu_{\u(n)}^{-1}(i\tilde\lambda)$. We are going to prove that in the dense open subset where $1/2 \pm \tilde\lambda > 0$ 
any $\Omega\in\mu_{\u(n)}^{-1}(i\tilde\lambda)$ is of the form $\Omega=k\Omega_0 k^t$, leading to the multiplicity freeness. From the restriction 
on $\tilde \lambda$, the matrix $A A^\dagger$ is invertible and we can write unambiguously the polar decomposition
$$
A = A_0 U_A,~ B= B_0 U_B, ~\Omega = -i A_0 U_A U_B^t B_0
$$
If we insert this decomposition in the relations of (\ref{parametrization_spn}) we get that
$$ U_A U_B^t \in U(n)_{i\tilde\lambda} ,~ U_A U_B^t = (U_A U_B^t)^t \;.$$

We want to show that we can find $k\in U(n)_{i\tilde\lambda}$ such that $U_A U_B^t= k k^t$ so that 
$\Omega = A_0U_A U_B^tB_0=A_0kk^tB_0=kA_0B_0k^t=k\Omega_0 k^t$. Indeed,  $U_A U_B^t$ can be diagonalized as
$U_A U_B^t = V u_0 V^\dagger$ with $u_0$ diagonal unitary matrix and $V$ unitary. By ordering the eigenvalues of $\tilde\lambda$,
$U_AU_B^t$ and so $V$ are block diagonal; in particular $V$ commutes with $\tilde\lambda$.
Since $U_AU_B^t$ is symmetric $V^t V$ commutes with $u_0$. It is always possible to choose a unitary 
square root $\sqrt{u_0}$ commuting
with $V^tV$. Then it is easy to check that $k=V \sqrt{u_0} V^\dagger$ is such that $U_A U_B^t= k^2$ with $k=k^t$ and $k\in U(n)_{i\tilde\lambda}$. \qed

\medskip
\begin{example}
{\rm Let us describe more explicitly the bihamiltonian polytope and the singularity locus in low dimension. We recall that the eigenvalues
are globally continuous functions and their derivative becomes singular on the border of the Weyl chamber of each subalgebra $\g_k$ appearing in (\ref{GC_spn}).

If $M_\phi=Sp(1)/U(1)$ then
${\cal C}(N_\phi)=\{\lambda\in\R|\, 0\leq\lambda\leq 2\}\sim \Delta_1$, the one dimensional simplex; $\lambda$ defines a smooth function. If $M_\phi=Sp(2)/U(2)$ then
${\cal C}(N_\phi)=\{(\lambda^{(k)}_i)_{1\leq i\leq 3-k\leq 2}\in\R^3|\, 0\leq\lambda^{(1)}_1\leq\lambda^{(2)}_1\leq\lambda^{(1)}_2\leq  2\}\sim \Delta_3$, the
three dimensional simplex. The Nijenhuis
eigenvalues $\lambda^{(1)}_i$ are singular when they reach the boundary of the positive Weyl chamber of $\g_1=\u(2)$ that happens when $\lambda^{(1)}_1=\lambda^{(1)}_2$.

Let $M_\phi=Sp(3)/U(3)$; then
$$
{\cal C}(N_\phi) = \{(\lambda^{(k)}_i) \in \R^6|\, 0\leq \lambda^{(1)}_1\leq\lambda^{(2)}_1\leq\lambda^{(1)}_2\leq\lambda^{(2)}_2\leq\lambda^{(1)}_3\leq 2,\,
\lambda^{(2)}_1\leq\lambda^{(3)}_1\leq\lambda^{(2)}_2\}\,.
$$
The singularity locus is reached on the boundary of the Weyl chamber of $\g_1=\u(3)$ and $\g_2=\u(2)$, that is when $\lambda^{(1)}_1=\lambda^{(1)}_2$,
$\lambda^{(1)}_2=\lambda^{(1)}_3$ or $\lambda^{(2)}_1=\lambda^{(2)}_2$. \qed}
\end{example}

\bigskip
\bigskip

\section{$M_\phi=SO(2n)/U(n)$}\label{sec_so2n}
Let us consider $M_\phi= SO(2n)/U(n)$. We showed in Section
\ref{sec_Rtacos} that the moment map in the fundamental representation of $\g_k=\u(n+1-k)\subset \so(2n)$ appearing in the chain 
(\ref{GC_so2n}) solves the master equation (\ref{master_equation}). By applying Proposition \ref{reduction_of_momentum_map} we get
the Njienhuis eigenvalues.

In the following theorem we prove that these are all the eigenvalues and that they are independent, by proving the complete integrability
of the collective hamiltonians defined by the chain (\ref{GC_so2n}).

\begin{theorem}\label{main_theorem_so_u}
The collective hamiltonians $F(\u(n),\ldots,\u(1))$ define a completely integrable model on $M_\phi=SO(2n)/U(n)$. The Nijenhuis tensor $N_\phi$
(\ref{Nijenhuis_tensor}) is of maximal rank and its ordered eigenvalues are
\begin{equation}
 \label{nijenhuis_eigenvalues_so_u}
 \lambda^{(k)}_i = 1 -2 \tilde\lambda^{(k)}_i \,, \;\;\; 1\leq i \leq n+1-k \leq n
\end{equation}
where $i\tilde\lambda^{(k)}_i$ are the eigenvalues of the moment map $\mu_{\u(n+1-k)}$. The bihamiltonian polytope is
the $n(n-1)/2$-dimensional ${\cal C}(N_\phi)\subset\R^{n(n+1)/2}$ where $(\lambda^{(k)}_i)\in{\cal C}(N_\phi)$ if, for $n=2N$,
$$
-1 \leq \lambda^{(1)}_{2N}=\lambda^{(1)}_{2N-1}\leq  \ldots \leq \lambda^{(1)}_2=\lambda^{(1)}_1 \leq 3,\qquad\quad
\lambda^{(k)}_i\leq\lambda^{(k+1)}_{i}\leq \lambda^{(k)}_{i+1}\;,
$$
and, for $n=2N+1$,
$$
-1 =\lambda^{(1)}_{2N+1}\leq \lambda^{(1)}_{2N}=\lambda^{(1)}_{2N-1} \leq  \ldots \leq \lambda^{(1)}_2=\lambda^{(1)}_1 \leq 3,\qquad\quad
\lambda^{(k)}_i\leq\lambda^{(k+1)}_{i}\leq \lambda^{(k)}_{i+1}\;.
$$
\end{theorem}
{\it Proof}.
If we write $X\in M_\phi$ in a block form with $X_{ij}\in M_n(\R)$ ,
$$
X= g \rho_\phi g^{-1} = \left(\begin{array}{cc}
                               X_{11} & X_{12} \cr
                               X_{21} & X_{22}
                              \end{array}\right)\,,\quad g\in SO(2n)\;,
$$
the moment map for the $\u(n)$ action is
$$
\mu_{\u(n)}(X) = X_{11}+X_{22}+ i (X_{12}-X_{21})\,,
$$
so that if $\tilde\lambda={\rm diag}(\tilde\lambda_1,\ldots,\tilde\lambda_n)$, a generic matrix in $\mu_{\u(n)}^{-1}(i\tilde\lambda)$ can be
written as
\bea  Z=
\left(\begin{array}{cc}
                               A & B+\tilde\lambda/2 \cr
                               B-\tilde\lambda/2 & -A
                              \end{array}\right),
\label{temp_V}\eea
with $A,B$ antisymmetric satisfying
\begin{equation}
[A,B]=[A,\tilde\lambda]=[B,\tilde\lambda]=0\,,\quad\quad A^2+B^2=\frac{1}{4}(\tilde\lambda^2-1)\,.\label{relABl}
\end{equation}
as a consequence of $X^2=\rho_\phi^2=-1/4$.

By using the Weyl group of $U(n)$ we can take $\tilde\lambda_{i-1}\leq \tilde\lambda_i$ so that 
$$\tilde{\lambda}=\textrm{diag}\,(\tilde \lambda_11_{m_1},\cdots ,\tilde \lambda_s1_{m_s})\;,$$ 
i.e. $\tilde\lambda_i$ has multiplicity $m_i$, $\sum m_i=n$. 
The condition (\ref{relABl}) implies that $A,\,B$ are block diagonal $A=\textrm{diag}\,(A_1,\cdots,A_s)$, $A_i\in M_{m_i}(\R)$, and the same goes for $B$. 
If $m_i$ is odd, then $A_i^2+B_i^2=1/4(\tilde\lambda_i^2-1)1_{m_i}$ implies $\tilde\lambda_i=\pm1$. We exclude first the possibility of $\tilde\lambda_i=-1$. 
Let $A_t=tA,\,B_t=tB$ and $\tilde \lambda_{it}=\textrm{sgn}\,\tilde\lambda_i\cdotp(1-t^2+t^2\tilde\lambda_i^2)^{1/2}$ (take $\textrm{sgn}\,(0)=\pm1$ does not matter), 
then $Z_t$ as given in (\ref{temp_V}) is a family of $\sp(n)$ matrices. Since $M_\phi$ is the orbit of $SO(2n)$, then the Pfaffian $pf(Z_t)$ is equal to
$pf(\rho_\phi)=(1/2)^n$; by evaluating it in $t=0$ we get
\bea \textrm{sgn}\prod\tilde\lambda_i^{m_i}=+1.\nn\eea
Hence $\tilde\lambda_i=-1$ with $m_i$ odd is excluded.

Thus one must always have even $m_i$, except possibly the last $m_s$ odd when $\tilde\lambda_s=1$. Thus for $n=2N$ all $m_i$ are even
while for $n=2N+1$, all but the last $m_i$ are even. We have then the Gelfand-Tsetlin pattern
\bea
-1\leq\tilde\lambda_1=\tilde\lambda_2\leq\cdots\leq\tilde\lambda_{n-1}=\tilde\lambda_n\leq1,~~n=2N\nn\\
-1\leq\tilde\lambda_1=\tilde\lambda_2\leq\cdots\leq\tilde\lambda_{n-2}=\tilde\lambda_{n-1}\leq\tilde\lambda_n=1,~~n=2N+1\nn\eea
To solve for $A,\,B$, we focus on the dense open subset where maximal amount of eigenvalues are distinct, thus $m_i=2$ for $n=2N$ case, and $m_i=2,~i=1,\cdots,N$, 
$m_{N+1}=1$ for $n=2N+1$. Then for both cases, one writes $A_i=a_i\sigma$, $B_i=b_i\sigma$, with $\sigma$ denoting the $2\times 2$ antisymmetric matrix, and
\bea a_i^2+b^2_i=\frac14(1-\tilde\lambda_i^2),~~i=1,\cdots, N.\label{temp_VI}\eea
To prove the complete integrability we use again Proposition \ref{Thimm_method} showing that for almost all $\tilde\lambda$ the action of the
stability subgroup $U(n)_{i\tilde\lambda}$ on $\mu_{\u(n)}^{-1}(i\tilde\lambda)$ is transitive.
From (\ref{temp_VI}), the orbits corresponding to fixed $\tilde\lambda$ are a product of $N$-circles. It is a direct check that the action
of $U(n)_{i\tilde\lambda}$
rotates $a_i+ib_i\to e^{2i\theta}(a_i+ib_i)$. The action is clearly transitive.
The same logic applies to $n$ odd. Finally to get (\ref{nijenhuis_eigenvalues_so_u}) one applies Proposition \ref{reduction_of_momentum_map}. \qed

\smallskip
\begin{rem}{\rm It is now straightforward to check that the number of independent eigenvalues described in the above theorem is the correct one. 
In the even case, e.g. for $n=4$, after renaming the independent eigenvalues the above pattern gives
\bea &&\left[\begin{array}{ccccccc}
      x_1 & ~ & x_1 & ~ & x_2 & ~ & x_2 \\
      ~ & x_1 & ~ & x_3 & ~ & x_2 \\
      ~ & ~ & x_4 & ~ & x_5\\
      ~ & ~ & ~ & x_6 & ~ \\ \end{array}\right],\nn\eea
and $\dim SO(8)/U(4)=28-16=12$. The counting for the general even case $n=2N$, when $\dim M_\phi= 2 N(2N-1)$, goes as
$$ N+(N-1)+\big(2N-2~+2N-3~+\cdots +1\big)=N(2N-1)\;.
$$
In the odd case, we have for $n=5$
\bea &&\left[\begin{array}{ccccccccc}
      1 & ~ & x_1 & ~ & x_1 & ~ & x_2 & ~ & x_2 \\
      ~ & x_3 & ~ & x_1 & ~ & x_4 & ~ & x_2 \\
      ~ & ~ & x_5 & ~ & x_6 & ~ & x_7\\
      ~ & ~ & ~ & x_8 & ~ & x_9 \\
      ~ & ~ & ~ & ~ & x_{10}  \\  \end{array}\right],\nn\eea
and $\dim SO(10)/U(5)=45-25=20$. The general counting for $n=2N+1$, when $\dim M_\phi =2 N(2N+1)$, goes as
\bea N+N+\big(2N-1~+2N-2~+\cdots +1\big)=N(2N+1)\;.\nn\eea
} 
\end{rem}

\section{$M_\phi=SO(n+2)/SO(n)\times SO(2)$}\label{sec_soodd}
Let us consider $M_\phi=SO(n+2)/SO(n)\times SO(2)$; in Section \ref{sec_Rtacos} we showed that the moment map of the subalgebra 
$$\g_k=\so(n+2-2k)\oplus\underbrace{\so(2)\oplus\ldots\so(2)}_k$$ 
in the representation $(S,i/2,\ldots,i/2)$ solves the master equation (\ref{master_equation}) so that every eigenvalue
defines a Nijenhuis
eigenvalue by (\ref{eigenvalues_moment_map_versus_nijenhuis}). We show in this section that by varying $k$ 
we get all Nijenhuis eigenvalues and that they are independent. 

We again make contact with the collective 
hamiltonians defined by (\ref{chain_subgroups_so}). This is equivalently described as the space of collective hamiltonians of the reduced chain 
\begin{equation}\label{chain_subgroups_so_reduced}
\so(2N+1) \supset \so(2N-1)\oplus\so(2) \supset \ldots \supset\so(3)\oplus\so(2) \supset \so(2) \supset 0,
\end{equation}
$$
\so(2N) \supset \so(2N-2)\oplus\so(2) \supset \ldots \supset\so(4)\oplus\so(2)\supset \so(2)\oplus\so(2) \supset 0~,
$$
where the $k$-th subalgebra of the chain $\g_{k}'=\so(n+2-2k)\oplus\so(2)$ is the subalgebra of $\so(n +2-2(k-1))\subset\g_{k-1}'$ corresponding to the non compact root $\alpha_1$. 
If $n+2=2N+1$ the last step 
is then $\g_{N}'=\so(2)$, if $n+2=2N$ the last step is then $\g_{N-1}'=\so(2)\oplus\so(2)$. We stress the fact that the difference between the chain
(\ref{chain_subgroups_so}) and (\ref{chain_subgroups_so_reduced}) is relevant only for the determination of the Nijenhuis eigenvalues and not 
for the definition of the collective hamiltonians.

\begin{thm}\label{main_theorem_so_so}
The collective hamiltonians $F(\so(n)\oplus\so(2),\so(n-2)\oplus\so(2),\ldots)$ define a completely integrable model on $M_\phi=SO(n+2)/SO(n)\times SO(2)$. Let
$n+2=2N$ or $2N+1$. The Nijenhuis tensor (\ref{Nijenhuis_tensor}) is of maximal rank and its eigenvalues are
\begin{equation}\label{bihamiltonian_hyperplane_so_odd}
\lambda^{(k)}_\pm = \pm |a_k| - \sum_{j=1}^k b_j + 1,\ k=1,\ldots N-1, \; {\rm and}\;\; \lambda^{(N)} = 1- \sum_{j=1}^Nb_j\ {\rm if}\;\; n+2=2N+1
\end{equation}
where $\pm i a_k$ are the eigenvalues of the moment map $\mu_{\so(n+2-2k)}$ for $\so(n+2-2k)\subset\g_k'$ and
$b_k =pf(\mu_{\so(2)})$
with $\so(2)\subset\g_{k}'$. 

The bihamiltonian polytope is then described as ${\cal C}(N_\phi)\subset\R^n$, where $(a_k,b_k)\in{\cal C}(N_\phi)$ if
\begin{equation}\label{ab_cone_so_odd}
0\leq |a_k|\leq |a_{k-1}|,~~~~ |b_k| \leq |a_{k-1}|-|a_k|~,  \;\;\;\; k=1,\ldots,N~,
\end{equation}
$a_0=1, a_N=0$ and, if $n+2 = 2N$, $b_N=0$ .
\end{thm}

{\it Proof}. Even though the Nijenhuis eigenvalues must be computed from the spin representations,
integrability of collective hamiltonians will depend on the properties of the moment maps of (\ref{chain_subgroups_so_reduced}) in the fundamental representation.
We are going first to characterize the coadjoint orbits contained in the image of the
moment map of the subalgebras appearing in (\ref{chain_subgroups_so_reduced}).
Let us parametrize $g\in SO(n+2)$ as
\bea g=\left(\begin{array}{ccc}
          \cdots & \vec \xi & \vec \eta \\
          \cdots & \vec x & \vec y \\
        \end{array}\right),~~\vec \xi,\vec \eta\in M_{n,1}(\R),~~\vec x,\vec y\in M_{2,1}(\R).\label{group_element}\eea
Since $\rho$ is written in block diagonal form as ${\rm diag}(0_{n},\sigma)$ we get
$$ \mu=g\rho g^{t} = \left(\begin{array}{cc}
          (\vec \xi,\vec \eta)\sigma(\vec \xi,\vec \eta)^t & (\vec \xi,\vec \eta)\sigma(\vec x,\vec y)^t \cr
          (\vec x,\vec y)\sigma (\vec \xi,\vec \eta)^t & (\vec x,\vec y)\sigma(\vec x,\vec y)^t
        \end{array}\right)=\left(\begin{array}{cc}
                                  h\, {\rm diag} (0_{n-2},a \sigma)h^t &  A\cr
                                  -A^t & b\sigma
                                 \end{array}\right) $$
where $b=\det (\vec x, \vec y)$, $a\in\R$ and $h\in SO(n)$. The last equality is just the standard form of a rank 2 antisymmetric matrix.
The reduction to $\h_\phi$ removes the off-diagonal blocks.

The $SO(n)\times SO(2)$ orbits contained in $\mu_{\h_\phi}(M_\phi)$ are
then the orbits ${\cal O}_{ab}$ through $\alpha_{ab}\equiv {\rm diag}(0_{n-2},a\sigma,b\sigma)$ parametrized by $a,b$.
Let $\h_{\alpha_{ab}}\subset\h_\phi$ the stability algebra of $\alpha_{ab}$.
If $a\not = 0$ then ${\cal O}_{ab}$ is isomorphic to the compact hermitian symmetric space of $SO(n)$ and
$h_{\alpha_{ab}}=\so(n-2)\oplus\so(2)\oplus\so(2)$. A generic point in $\mu_{\h_\phi}^{-1}(\alpha_{ab})$ is of the form
\bea
P=\left( \begin{array}{ccc} 0_{n-2}&0&0\cr
                          0 & a\sigma & X\cr
                          0&-X^t & b\sigma \end{array}\right)\label{temp_VII}\eea
where $X = (\vec u,\vec v)\sigma(\vec x,\vec y)^t\in M_2(\R)$, with $\vec u,\vec v,\vec x,\vec y\in \R^2$ satisfying  
$\det (\vec u,\vec v)=a$, $\det (\vec x,\vec y)=b$ and
$$
g_4 = \left(\begin{array}{cc}
           \vec u & \vec v \cr
           \vec x & \vec y
          \end{array}
 \right) \;\;\;  g_4 g^t_4=1\;.
$$

The action of $(g_{n-2},h,k)\in SO(n-2)\times SO(2)\times SO(2)$ integrating $\h_{\alpha_{ab}}$ is given
by $X\rightarrow h X k^t$. By combining left $SO(2)$ action on $X$ and reparametrization of $(\vec{x},\vec{y})$
we can choose $(\vec u,\vec v)=\textrm{diag}\,(u,v)$; indeed we can choose $h,k\in SO(2)$ such that
$$
X = (\vec{u},\vec{v})\sigma (\vec{x},\vec{y})^t = h h^t(\vec{u},\vec{v})k \sigma k^t (\vec{x},\vec{y})^t =
h\, \textrm{diag}\,(u,v)\sigma(\vec{p},\vec{q})^t\;.
$$

Orthogonality of $g_4$ then means
\bea u^2+|\vec p|^2=1=v^2+|\vec q|^2,~~\vec p^t\vec q=0.\nn\eea
Let $\vec p=(p_1,p_2)\not=0$, then $\vec q=c(-p_2,p_1)$ for some $c$. Since $b=\det(\vec p,\vec q)=c|\vec p|^2$ and $a=uv$, we get
\bea\label{equation_ab} 1+a^2-b^2= u^2+\frac{a^2}{u^2}.\eea
The condition that there are real solutions for $u^2$, together with the upper bound of $|a|$, gives the range of $(a,b)$
\bea \label{range_ab} b^2\leq (1-|a|)^2,~~~~ |a|\leq 1.\eea
The space of solutions to the above equation is just the space of those $\vec p\in\R^2$ with $|\vec p|^2=1-u^2$ and the right $SO(2)$ action on $X$ is transitive on this
circle. We conclude that the action of $\h_{\alpha_{ab}}$ on $\mu_{\h_\phi}^{-1}(\alpha_{ab})$
is transitive.

If $a=0$ then $\h_{\alpha_{0b}}=\h_\phi$; moreover $\vec{\xi},\vec{\eta}$ appearing in (\ref{group_element})
are collinear and it can be shown that $|b|\leq 1$, extending (\ref{range_ab}) to the case $a=0$.

The orbits of the subgroups appearing in the two chains of
(\ref{chain_subgroups_so_reduced}) will have the same pattern, compact hermitian symmetric spaces or points. We get two new variables
$(a_k,b_k)$ for each step, until we get
to $\g_{N-1}'$, which is the last step for the even case. In the odd case there is one more reduction to $\g_{N}'=\so(2)$ that gives us one more $b_{N}$ variable.
In both cases we get
\bea n=\frac12\dim SO(n+2)-\frac12\dim (SO(n)\times SO(2))\nn\eea
variables, which is consistent.
In order to establish the range of these variables let ${\cal O}_{a_kb_k}\subset p_{k-1,k}({\cal O}_{a_{k-1}b_{k-1}})$, 
where $p_{k-1,k}:\g_{k-1}'\rightarrow \g_{k}'$ is the dual of the inclusion map, denote the adjoint orbit of 
$\alpha_{a_k,b_k}={\rm diag}(0_{n-2k},a_k\sigma,b_k\sigma)\in\g_k'$. Then, since ${\cal O}_{a_{k-1}b_{k-1}}$ is isomorphic to the $SO(n+2-2(k-1))$ orbit
of $a_{k-1}\rho_{k-1}$, where $\rho_{k-1}$ is the normalized generator of the non compact root, we repeat the above considerations and conclude that
$(a=a_k/a_{k-1},b=b_k/a_{k-1})$ satisfy inequalities (\ref{range_ab}) and
so (\ref{ab_cone_so_odd}). Moreover, we showed above that the action of the stability group of $\alpha_{a_k,b_k}\in\g_k'$ is transitive on 
$p_{k-1,k}^{-1}(\alpha_{a_k,b_k})$. 
By applying Proposition \ref{Thimm_method} we prove the complete integrability.

Finally we have to compute the eigenvalues of the moment map of $\g_k$ in the representation $(S,i/2,\ldots,i/2)$ in terms of $a_k,b_k$. 
Since the weights of the spin representation are $(\pm 1/2,\ldots,\pm 1/2)$ they are easily
computed as
\bea \pm\frac{i}{2}a_k+\frac{i}{2}\sum_{j=1}^kb_j\;.\nn\eea
By using Proposition \ref{reduction_of_momentum_map}, these lead to the pointwise eigenvalues of $N_\phi^*$
\bea
\pm a_i - \sum_{j=1}^ib_j + 1~~~~~ i=1,\ldots n. \qed\nn\eea

\begin{rem}{\rm
Note that the spin representation for $n$ even is reducible, but it does not have any effect on the proof. Also for $n$ even, the last reduction  
$\so(4)^*\to\so(2)^*\oplus\so(2)^*$ does not take place through removing the root $\alpha_1$ as the earlier steps, but this again has no effect on the validity of the proof. }
\end{rem}

\smallskip
\begin{rem}{\rm
Let us identify the hamiltonians of the action of the Cartan subalgebra $\t\subset\so(n+2)$. Indeed, the $b_k$ are the hamiltonians of the $(n+2-k)$-th $\so(2)$,
$k\leq N$ when $n+2=2N+1$ and $k\leq N-1$ for $n+2=2N$. In the even case, the missing generator is given by the last $a_{N-1}$. These variables are of course 
global smooth functions.}
\end{rem}
\smallskip
\begin{rem}{\rm 
The value of $a$ appearing in (\ref{temp_VII})
can always be assumed to be non-negative, except in the even case in the last step $\so(4)^*\to\so(2)^*\oplus\so(2)^*$.
Indeed, conjugating $P$ by a rotation of $\pi$  along, say, the $(n-2)$, $(n-1)$ direction flips $a\to -a$. If we think to the definition 
of the action variables described at the end of Section \ref{collective_complete_integrability},  
then $|a_k|$ is obtained by projecting $\mu_{\so(n+2-2k)}$ to the positive Weyl chamber. 

In the last even step it is then convenient not to introduce the absolute value in the definition of the Njenhuis eigenvalue, since $a_{N-1}$ and then $\lambda_\pm^{(N-1)}$ are
smooth global functions while the absolute value would introduce a singularity.

}
\end{rem}

\section{Conclusions}

In this paper we proved that the PN structures defined on compact hermitian symmetric spaces are of maximal rank, or equivalently that they define a completely
integrable model that admits a bihamiltonian description. In the case of Grassmannians we recover the well known Gelfand-Tsetlin integrable model, so that our result
can be phrased by saying that we show that Gelfand-Tsetlin variables are in involution also with respect to the Bruhat-Poisson structure. In the other cases,
the results are new also on the symplectic side. From our point of view, it is natural to look for the information about the Poisson pencil that are
contained in these models. We collect here some observations that we plan to develop in the future.

\medskip

\noindent 1) {\it Geometry of the Poisson pencil and log symplectic structures}. The description of the spectrum of the Nijenhuis tensor $N_\phi$ gives information on
the geometry of the pencil $\pi_t = \pi_0+t\Omega^{-1}$, where $\pi_0$ is the Bruhat-Poisson structures and $\Omega$ is the KKS symplectic form. We collect here
few basic observations.

The knowledge of eigenvalues allows to reconstruct the strata of symplectic leaves of a given dimension. In fact, the corank of $\pi_t$ at a given point is the multiplicity
of the eigenvalue $-t$ so that the symplectic foliation can be analyzed by means of the hyperplanes ${\cal C}^{(k)}(t)$ of ${\cal C}(N_\phi)$ defined as the set of points where
the $k-th$ eigenvalue is equal to $-t$. For instance we can conclude that
on the complement of the preimage of ${\cal C}(t)=\bigcup_k {\cal C}^{(k)}(t)$ $\pi_t$ is nondegenerate; in particular $\pi_t$ is the inverse of a symplectic form for
all $t$ bigger than the radius
of the smallest ball containing ${\cal C}(N_\phi)$. This behaviour is a clear hint of a {\it log symplectic structure}, that we plan to discuss in a separate paper. 
In particular, we plan to investigate the relation with the framework
of {\it tropical moment map} introduced in \cite{Gualtieri} and the very recent \cite{KMS}.

Moreover, as described at
the end of Subsection \ref{Poisson-Nijenhuis structures}, for each $t$ the modular vector field of $\pi_t$ with respect to the symplectic volume form is given by
the symplectic vector field $\Omega^{-1}_{kks}d{\rm Tr}N_\phi$. This vector field
is not hamiltonian in general for $\pi_t$, but it is easy to see that $\log\det (N_\phi+t)$ gives a local hamiltonian, which is well defined provided that no
Nijenhuis eigenvalue is equal to $-t$.

\medskip

\noindent 2) {\it Lifting to the symplectic groupoid}. In \cite{BCQT} the Poisson Nijenhuis structure on $\C P_n$ was used to quantize the symplectic groupoid of
the Bruhat-Poisson structure. As briefly summarized in the Introduction, the procedure requires the integration to a groupoid cocycle of the Poisson vector field
$\Omega^{-1}_{kks}d\lambda$ associated to every Nienhuijs eigenvalue $\lambda$. This gives an integrable model on the symplectic groupoid compatible with the multiplication.
In this construction, it is crucial that the eigenvalues are smooth global functions. In general we know that the Nijenhuis eigenvalues are
globally continuous functions but their differential becomes singular on the boundary of the Weyl chamber of each of the subalgebras appearing in
(\ref{chain_subalgebras_general}). So the singularity locus can be read from our construction and this analysis will be done in a separate paper.
In general, it is an interesting problem to put this peculiar procedure of integration
of cocycles under the light of the more canonical integration of Poisson Nijenhuis structures developed in \cite{SX}.

\end{document}